\pdfoutput=1
\documentclass[onefignum,onetabnum]
{siamart251216}

\usepackage{amsmath,amssymb}
\usepackage{graphicx}
\usepackage{booktabs}
\usepackage{array}
\usepackage{mdframed}
\usepackage{algorithm}

\newsiamremark{remark}{Remark}
\newsiamthm{conjecture}{Conjecture}

\crefname{definition}{Definition}{Definitions}\Crefname{definition}{Definition}{Definitions}
\crefname{proposition}{Proposition}{Propositions}\Crefname{proposition}{Proposition}{Propositions}
\crefname{lemma}{Lemma}{Lemmas}\Crefname{lemma}{Lemma}{Lemmas}
\crefname{corollary}{Corollary}{Corollaries}\Crefname{corollary}{Corollary}{Corollaries}
\crefname{remark}{Remark}{Remarks}\Crefname{remark}{Remark}{Remarks}
\crefname{conjecture}{Conjecture}{Conjectures}\Crefname{conjecture}{Conjecture}{Conjectures}
\AddToHook{env/definition/begin}{\crefalias{theorem}{definition}}
\AddToHook{env/proposition/begin}{\crefalias{theorem}{proposition}}
\AddToHook{env/lemma/begin}{\crefalias{theorem}{lemma}}
\AddToHook{env/corollary/begin}{\crefalias{theorem}{corollary}}
\AddToHook{env/remark/begin}{\crefalias{theorem}{remark}}
\AddToHook{env/conjecture/begin}{\crefalias{theorem}{conjecture}}

\newcommand{\dname}[1]{{\upshape(#1.)}\ }

\newcommand{\Uex}{\textsf{U-Example}}
\newcommand{\Uexs}{\textsf{U-Examples}}
\newcommand{\eps}{\varepsilon}
\newcommand{\Ng}{N_{\mathrm{g}}}
\newcommand{\Nr}{N_{\mathrm{r}}}
\newcommand{\Nz}{N_{0}}
\newcommand{\Sstar}{S^{\star}}
\newcommand{\R}{\mathbb{R}}
\newcommand{\Bl}{\mathrm{B}}
\newmdenv[linecolor=black!25,linewidth=0.6pt,roundcorner=3pt,
  backgroundcolor=black!2,innertopmargin=6pt,innerbottommargin=6pt,skipabove=0pt,skipbelow=0pt]{payloadbox}
\newenvironment{payload}{\par\addvspace{7pt}\begin{payloadbox}}{\end{payloadbox}\par\addvspace{7pt}}

\newcommand{\TheTitle}{U-Examples: A Unified Family of Worst-Case Instances for Optimization Algorithms}

\title{{\TheTitle}}
\author{Eghbal Hosseini\thanks{Department of Applied Mathematics and Computer Science, Technical
University of Denmark, 2800 Kongens Lyngby, Denmark (\email{eghho@dtu.dk}).}}
\headers{U-Examples: worst-case instances for optimization algorithms}{E. Hosseini}

\begin{document}
\maketitle

\begin{abstract}
Worst-case instances show where optimization algorithms fail: the Klee--Minty cube for the simplex
method, ill-conditioned and rippled bowls for gradient descent, degenerate minima for Newton's
method, redundant Klee--Minty cubes for interior-point methods, unprunable trees for
branch-and-bound, deceptive landscapes for metaheuristics, and degenerate cycling for active-set
methods. These constructions were developed over seventy years in separate literatures, each with
its own tools. This survey reviews them together. We trace the history of rule-specific worst cases
across seven representative algorithm classes and organize them with one reading: in each, the optimum lies
near the start, yet the algorithm's own update rule is misdirected onto a long trajectory or, at the
boundary, into a non-terminating loop. We call such an instance a \emph{\Uex}; the letter U stands
both for the unification and for the U-shaped bowls around which descent methods are made to crawl.
A common framework treats every algorithm as a rule on a state space with a progress potential, so
that one definition of \emph{misdirection}, and of the resulting \emph{blow-up} in cost, covers
pivoting, descent, interior-point, tree-search, and population methods. Each member is drawn in a
common visual grammar and quantified in general dimension $n$. The blow-up is independent of the
dimension for the two descent members and exponential in it for four search members, with
degenerate cycling as a breakdown. We connect the constructions to neighbouring theories---oracle
complexity, performance estimation, smoothed analysis, and the complexity of local search---each of
which asks a different question of the same instances. We then map each worst case to the remedy
that defeats it, examine how robust the instances are, and close with open challenges, in
particular whether a worst case can be predicted from the update rule alone.
\end{abstract}



\section{Introduction}\label{sec:intro}
Optimization is the problem of choosing, among all admissible decisions, one that is best: the
cheapest schedule, the most likely parameters of a model, the lightest structure that carries its
load. It is one of the main routes by which mathematics reaches applications in engineering,
economics, statistics, and machine learning, and the problems it poses differ widely in kind. The
variables may be continuous or discrete, the objective linear or nonlinear, smooth or rugged, the
constraints few or many. No single method serves them all. Over seventy years the field has
developed a family of algorithm classes, each fitted to one kind of problem: the simplex and
interior-point methods for linear programs, gradient and Newton methods for smooth unconstrained
problems, active-set methods for constrained nonlinear ones, branch-and-bound for integer programs,
and metaheuristics for problems with little usable structure. These seven classes, the core of the
standard curriculum, are the subject of this survey.

Different as they are, these methods share a form. Each is an iteration: from its current state it
produces the next one by a fixed local rule---follow the negative gradient, minimize a local
quadratic model, move to a better neighbouring vertex, split a problem in two, resample around the
best candidates found so far. The rule sees only local information, and the basic question about
any algorithm is how many applications of its rule a problem requires. There are three standard
ways to answer. \emph{Worst-case} analysis bounds the number of steps on the hardest input of a
given size and yields a guarantee. \emph{Average-case} analysis takes the expectation over a
distribution of inputs. \emph{Smoothed} analysis takes the worst case over inputs after each has
been slightly perturbed at random. The three can disagree sharply. The simplex method is the
best-known example: it solves linear programs efficiently in practice, although almost every pivot
rule that has been analysed needs exponentially many steps on some input (\cref{sec:history}).

Worst-case instances are the witnesses of this gap, and they matter more than their rarity
suggests. An upper bound says that a method never costs more than a stated amount; only an explicit
instance shows that the bound cannot be improved, and where the reasoning of the method fails. Each
instance exhibits a precise way in which a rule can be fooled, and much of algorithm design can be
read as a sequence of answers to them: anti-cycling rules, acceleration, regularization, cutting
planes, and restarts were each introduced to defeat a known bad case.

Such instances are hard to find. A method that works well on almost every input gives few clues
about the input on which it fails, and the known constructions are ingenious and specific.
Optimization is organized by algorithm class, and so is its folklore of hard instances. Each class
relies on some property of the problem to make progress, and its hard instance is built to break
exactly that property. The simplex method has the Klee--Minty cube, on which Dantzig's pivot rule
visits all $2^n$ vertices \cite{kleeminty1972}. Gradient descent has ill-conditioned quadratics, on
which the iterate zig-zags across a narrow valley \cite{nesterov2004}. Newton's method has
degenerate minima, at which its quadratic convergence drops to linear \cite{griewank1985}.
Interior-point methods have Klee--Minty cubes with redundant constraints, which bend the central
path past every vertex \cite{deza2008}. Branch-and-bound has integer programs whose relaxations
prune nothing \cite{jeroslow1974}. Evolutionary and other population-based metaheuristics have
deceptive landscapes that draw the search away from the optimum \cite{deb1993,droste2002}.
Active-set methods have degenerate points, at which the working set can change forever while the
iterate stands still \cite{beale1955}. These constructions arrived one at a time, some of them
decades after the method they defeat (\cref{sec:history}). Each lives in its own literature, with
its own notation, its own extremal construction, and its own intuition. There is no general recipe
for producing them, and to our knowledge they have not been surveyed together.

\begin{figure}[t]\centering
\includegraphics[width=\textwidth]{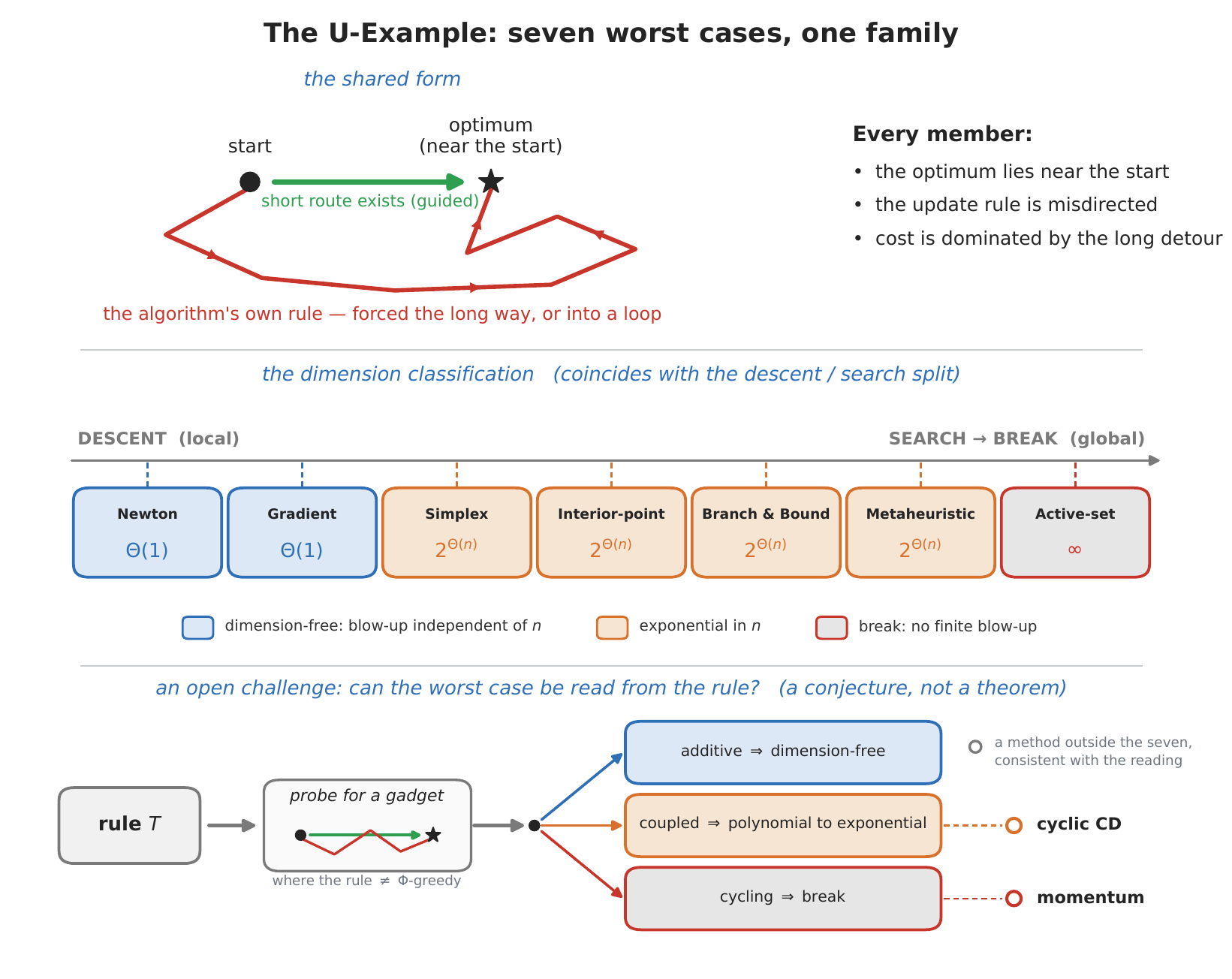}
\caption{Overview of the \Uex\ reading in three registers. \emph{Top:} the shared form---the
optimum lies near the start and a short guided route exists, yet the algorithm's own rule is forced
the long way, or, at the boundary, into a loop. \emph{Middle:} the dimension classification---the
seven worst cases along the descent/search axis, coloured by \cref{thm:classification}
(dimension-free for Newton and gradient; exponential for simplex, interior-point, branch-and-bound
and metaheuristic; the active-set break). \emph{Bottom:} the open challenge of
\cref{sec:outlook}---reading a worst case from the rule alone: a rule is probed for a gadget whose
composition (additive, coupled, or cycling) suggests its class (\cref{conj:general}); momentum and
cyclic coordinate descent are two methods outside the seven that are consistent with the reading.
The seven members are drawn out, one per row, in \cref{fig:plate}.}
\label{fig:story}
\end{figure}

The prevailing reading of each is that it encodes a \emph{hard problem}. We adopt a different
reading, and it is the organizing idea of this survey: each encodes a \emph{misdirected rule}. In
every case the optimum sits close to the start---a short route exists---but the algorithm's update
rule points, step after step, somewhere other than along that route. The Klee--Minty cube is not a
hard linear program: its optimum is a vertex adjacent to the starting vertex, and a different pivot
reaches it in one step. What is exponential is the behaviour of one pivot rule on that cube. A
flat-bottomed bowl is not a hard function to minimize: the minimizer is one correctly scaled step
away, and Newton's method simply mis-scales its step. We call an instance of this kind a
\emph{\Uex}. The reading does not remove the difficulty of constructing worst cases; whether they
can be predicted from the rule alone is the open challenge of \cref{sec:outlook}. What it does is
put the known ones in a common form, in which they can be defined once, measured in the same
units, compared across classes, and matched to the remedies that defeat them. \Cref{fig:story}
summarizes the reading in three registers: the form the worst cases share, how their cost scales
with the dimension, and the open question of predicting them.

\subsection{Scope and organization}
The survey covers \emph{rule-specific} worst cases: instances built against one update rule, as
opposed to lower bounds that hold for every method of a class. It treats seven representative
classes in detail and points to the corresponding results for other methods. It is organized as
follows. \Cref{sec:history} traces the history of such constructions across the strands of
optimization in which they arose. \Cref{sec:framework} sets up a common framework in which every
algorithm is a rule on a state space with a progress potential, so that one definition of
misdirection and of blow-up applies to pivoting, descent, interior-point, tree-search, and
population methods. \Cref{sec:members} treats the seven members in turn; for each it reviews the
literature, states the \Uex, computes its cost in general dimension~$n$, and summarizes the current
state of the art; it closes with the rule-specific worst cases of methods beyond the seven
(\cref{ssec:beyond}). \Cref{sec:theorem} classifies the members by how the blow-up scales with the
dimension. \Cref{sec:synthesis} compares the members across classes: what they share, and the axes
along which they differ. \Cref{sec:connections} relates the constructions to neighbouring
fields---smoothed analysis, oracle complexity, performance estimation, the complexity of local
search, and the optimizers of machine learning. \Cref{sec:remedy} maps each worst case to the
algorithmic remedy that defeats it, and \cref{sec:boundaries} discusses where the unifying reading
is exact, where it strains, and how robust the constructions are. \Cref{sec:outlook} collects open
challenges, chief among them the converse question of whether a worst case can be read off from the
update rule alone; we state what can be proved, pose a conjecture, and examine two methods outside
the seven. \Cref{sec:conclusion} concludes, and \cref{app:proofs} contains the proofs.

\subsection{Relation to prior work}
We are not aware of earlier work that organizes the worst-case instances of these classes under a
single reading. Two bodies of work come closest, and \cref{sec:connections}
discusses both. The performance estimation problem computes the exact worst case of a first-order
method on a class of convex functions by semidefinite programming
\cite{drori2014,taylor2017,taylor2017pesto,goujaud2022pepit}; it is computational and specific to
first-order methods, whereas the reading used here is geometric and spans discrete, combinatorial,
and stochastic methods. Unifications in the literature link two classes at a time through polytope
geometry: deformed products for the simplex worst cases \cite{amentaziegler1999,ziegler1995}, and
central-path curvature for the simplex and interior-point cubes \cite{deza2009curvature}. The
\Uex\ reading is of a different kind. It rests on the misdirection of the update rule rather than on
polytope structure, and it applies to all seven classes.

\subsection{Novelty and scope}
The survey contributes at three levels.

\emph{The reading and the framework.} The idea that a worst-case instance encodes a misdirected
rule, not a hard problem, is the starting point. The framework of \cref{sec:framework} makes it
precise for pivoting, descent, interior-point, tree-search, and population methods alike, and the
accounting bound of \cref{prop:blowup} names the two ways in which a rule can be forced the long
way. Both are introduced here.

\emph{The \Uex\ form of the seven members.} For every class we construct the same object: a normal
instance and its \Uex\ side by side, the control parameter that turns one into the other, and the
three costs, drawn in one visual grammar (\cref{fig:plate,fig:plate3d}) and computed in general
dimension (\cref{sec:members}). For four classes---Newton's method, gradient descent,
interior-point methods, and metaheuristics---a formulation of this kind has, to our knowledge, not
been given before. One of the four is also a new instance: the rippled bowl for gradient descent
(\cref{ssec:gradient}), introduced here with a proof of its cost.

\emph{The cross-class results.} The classification by dimension (\cref{thm:classification}), the
synthesis along four axes (\cref{sec:synthesis}), the map from each worst case to its remedy
(\cref{sec:remedy}), and the computed study of robustness (\cref{sec:boundaries}) compare the
members with one another, which the separate literatures could not do.

The instances and phenomena beneath these formulations are classical, with the one exception just
named, and we use them with attribution: the Klee--Minty cube for the simplex method
\cite{kleeminty1972}; its redundant version, due to Deza, Nematollahi, and Terlaky, for
interior-point methods \cite{deza2008}; Jeroslow's integer program for branch-and-bound
\cite{jeroslow1974}; degenerate cycling for active-set and pivoting methods \cite{beale1955}; the
slow convergence of Newton's method at a degenerate minimizer \cite{deckerkelley1980,griewank1985};
and the multimodal landscapes of global optimization, related to the deceptive landscapes of
\cite{deb1993}. The mechanism of the rippled bowl---the step length of the method is dictated by
the largest curvature anywhere---is textbook as well. In keeping with the purpose of a survey, the
new results serve to unify and illustrate these existing ideas. Where a member's cost is
probabilistic (metaheuristics) or its blow-up is not finite (active-set), we say so and adjust the
claim rather than force it into a uniform mould.

\Cref{tab:overview} previews the family: the seven worst-case instances, the control parameter that
switches each one on, the mechanism by which the rule is misdirected, and the resulting blow-up.

\begin{table}[t]
\renewcommand{\arraystretch}{1.35}\footnotesize\centering
\caption{The family at a glance. Seven worst-case instances, one form: the optimum is near, the
rule is misdirected, and the cost is dominated by the long detour. The \emph{control parameter} is
the feature of the instance that switches the worst case on. The last column anticipates the
dimension classification of \cref{thm:classification}.}
\label{tab:overview}
\begin{tabular}{@{}p{1.85cm} >{\raggedright\arraybackslash}p{2.05cm} >{\raggedright\arraybackslash}p{1.75cm} >{\raggedright\arraybackslash}p{2.65cm} >{\raggedright\arraybackslash}p{2.4cm}@{}}
\toprule
\textbf{Member} & \textbf{Worst-case instance} & \textbf{Control parameter} & \textbf{How the rule is misdirected} & \textbf{Blow-up (class)} \\
\midrule
Newton & flat-bottomed bowl $\sum_i x_i^{2m}$ & flatness $2m$ & the curvature model under-scales the step & $\approx\!(2m{-}1)\ln\tfrac1\eps$ \newline (dimension-free) \\
Gradient & rippled bowl & roughness $\rho$ & ripples dictate a tiny step length & $\Theta(\rho\ln\tfrac1\eps)$ \newline (dimension-free) \\
Simplex & Klee--Minty cube & dimension $n$ (and scaling) & the pivot rule is lured onto a $2^n$-vertex path & $2^n-1$ \newline (exponential) \\
Interior-point & redundant Klee--Minty cube & redundancy, dimension $n$ & redundant walls bend the central path & $2^{\Theta(n)}$ \newline (exponential) \\
Branch-and-bound & Jeroslow's integer program & parity of $n$ & the relaxation bound prunes nothing & $2^{\Theta(n)}$ \newline (exponential) \\
Metaheuristic & landscape with $k^n$ basins & basins per coordinate $k$ & basins split the population & $k^n$ expected \newline (exponential) \\
Active-set & degenerate point & degeneracy & the working set cycles, the iterate is fixed & $\infty$ \newline (break) \\
\bottomrule
\end{tabular}
\end{table}

\section{A brief history of worst-case analysis in optimization}\label{sec:history}
\begin{figure}[t]\centering
\includegraphics[width=\textwidth]{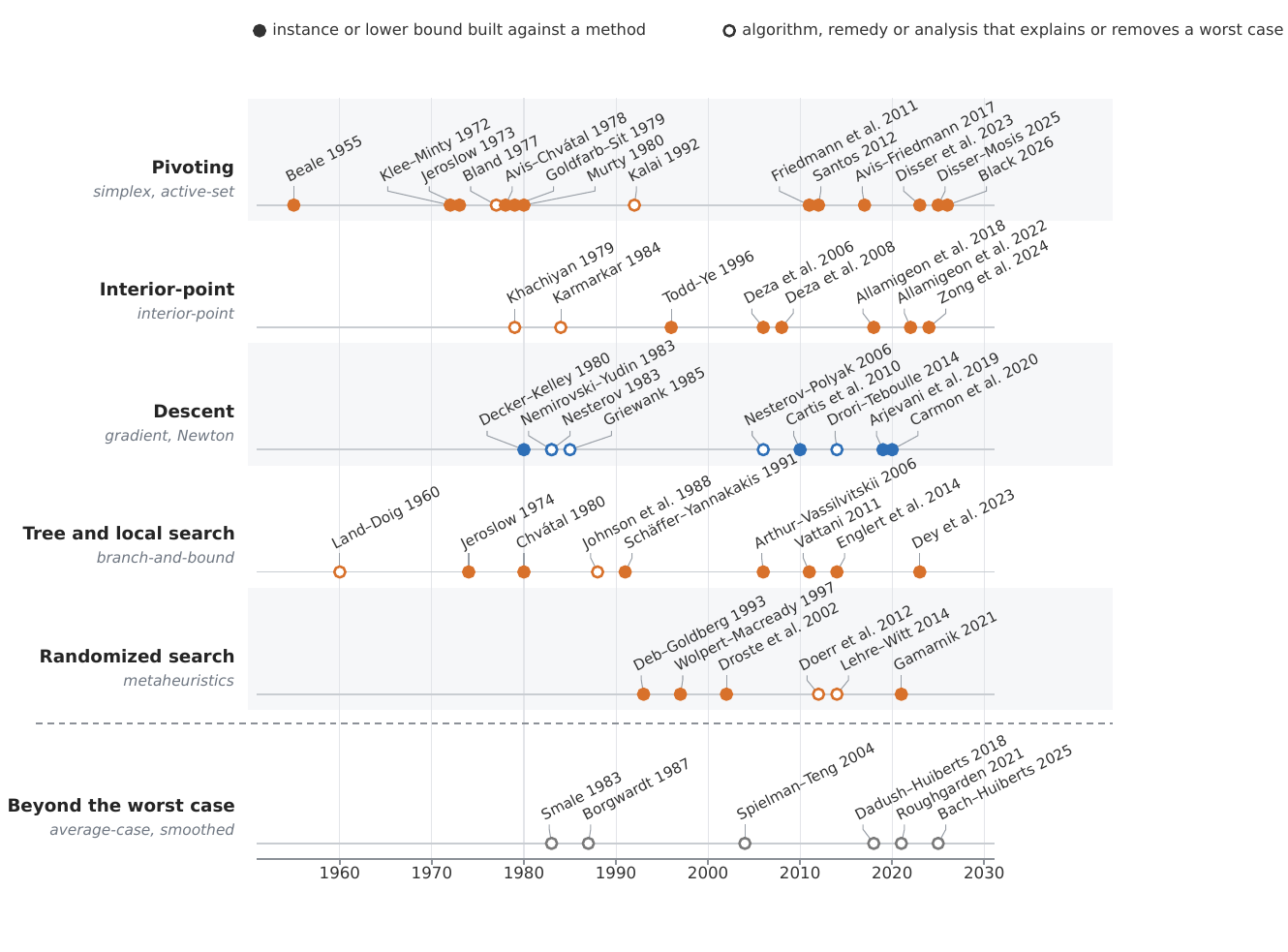}
\caption{Seventy years of worst-case constructions, by strand. Each of the upper five lanes is a
strand of \cref{tab:strands}, with the classes it serves named beneath it; a lane is coloured by
the dimension class of those classes (blue: dimension-free; orange: exponential;
\cref{thm:classification}), as in the other figures. Filled markers are instances or lower bounds
built against a method; open markers are algorithms, remedies, or analyses that explain or remove
a worst case. The bottom lane (grey) is not an algorithm class: it collects the average-case and
smoothed analyses of \cref{ssec:h-beyond}, almost all of which concern the simplex method. Every
entry is discussed, with its reference, in \cref{sec:history}.}
\label{fig:timeline}
\end{figure}

Worst-case analysis asks how much work an algorithm needs on the hardest input of a given size. It
became the organizing yardstick of algorithm analysis in the 1960s, when running time bounded by a
polynomial in the input size was proposed as the formal meaning of an efficient algorithm
\cite{edmonds1965}. In optimization the yardstick plays two different roles. As an \emph{upper bound}
it certifies a method: a guarantee that, on every instance, a stated number of iterations suffices.
As a \emph{lower bound} it exposes a method, and a lower bound for a specific method is almost always
witnessed by a construction---an explicit, usually small and highly structured instance on which the
method's rule is provably slow. This survey is about the second role. The constructions arose
independently in five strands of optimization. The strands are not our seven algorithm classes: a
strand may serve two classes, and one of them also covers local search, which lies outside the
seven. \Cref{tab:strands} gives the correspondence, and \cref{fig:timeline} places the strands on a
common time axis. This section follows them in turn (\cref{ssec:h-pivot,ssec:h-ipm,ssec:h-descent,ssec:h-tree,ssec:h-random})
and closes with the work that looks beyond the worst case (\cref{ssec:h-beyond}); the seven classes
themselves are treated one by one in \cref{sec:members}.

\begin{table}[t]
\renewcommand{\arraystretch}{1.2}\footnotesize\centering
\caption{The five historical strands of this section and the algorithm classes they serve.}
\label{tab:strands}
\begin{tabular}{@{}lll@{}}
\toprule
\textbf{Strand} & \textbf{Algorithm classes served} & \textbf{Treated as members in} \\
\midrule
Pivoting methods (\S\ref{ssec:h-pivot}) & simplex, active-set & \S\ref{ssec:simplex}, \S\ref{ssec:activeset} \\
Interior-point methods (\S\ref{ssec:h-ipm}) & interior-point & \S\ref{ssec:ipm} \\
Descent methods (\S\ref{ssec:h-descent}) & gradient, Newton & \S\ref{ssec:gradient}, \S\ref{ssec:newton} \\
Tree search and local search (\S\ref{ssec:h-tree}) & branch-and-bound & \S\ref{ssec:bb} \\
Randomized search (\S\ref{ssec:h-random}) & metaheuristics & \S\ref{ssec:meta} \\
\bottomrule
\end{tabular}
\end{table}

\subsection{Pivoting methods: simplex and active-set}\label{ssec:h-pivot}
The simplex method moves from vertex to adjacent vertex of the feasible polytope, choosing each
improving edge by a \emph{pivot rule} \cite{dantzig1963}. Its first worst case was not slowness but
non-termination: at a degenerate vertex, where more constraints are active than the dimension
requires, a pivot rule can cycle through bases forever without moving, as Beale's example shows
\cite{beale1955}. Bland's smallest-subscript rule restored termination for all instances
\cite{bland1977}. Termination is not efficiency, however. Klee and Minty constructed a slightly
deformed cube on which Dantzig's original rule visits all $2^n$ vertices \cite{kleeminty1972}, and
exponential examples followed for essentially every deterministic rule that was examined: the
greatest-improvement rule \cite{jeroslow1973}, Bland's rule itself \cite{avischvatal1978}, the
steepest-edge rule \cite{goldfarbsit1979}, and parametric linear programming, on which the
shadow-vertex rule rests \cite{murty1980}. The constructions share a template, the \emph{deformed
product} of \cite{amentaziegler1999}: a polytope that is combinatorially a cube, with facets tilted
so that a monotone path through all vertices exists and the rule is lured onto it.

The template has proved remarkably durable. History-dependent rules, designed to defeat it, fell
much later: Zadeh's least-entered rule after a subexponential lower bound \cite{friedmann2011zadeh}
and then an exponential one \cite{disser2023zadeh}, and Cunningham's rule
\cite{avisfriedmann2017}. Most recently, Black constructed linear programs on which the shadow-vertex
rule is exponential for every choice of projection, and the steepest-edge rule for every choice of
norm \cite{black2026}. Randomization helps but does not settle the matter: random-facet rules admit
subexponential upper bounds \cite{kalai1992,matousek1996}, subexponential lower bounds are also known
\cite{friedmannhansenzwick2011}, and on Klee--Minty cubes specifically the number of pivots of the natural
random rules grows only about quadratically with $n$ \cite{gartnerhenkziegler1998}. Whether \emph{any} pivot rule runs in
polynomial time is open. The closely related question of whether linear programming admits a
strongly polynomial algorithm is the ninth of Smale's problems for the twenty-first century
\cite{smale1998}, and the geometric question underneath---how long a shortest edge path between two
vertices of a polytope can be---survived the refutation of the Hirsch conjecture \cite{santos2012}
with only quasi-polynomial upper bounds \cite{kalaikleitman1992,todd2014}.

The simplex method is the linear special case of a larger class. \emph{Active-set} methods for
quadratic and nonlinear programming maintain a working set of constraints treated as equalities and
change it one constraint at a time \cite{nocedalwright2006,fletcher2000}; applied to a linear
objective they are the simplex method. They inherit both of its failures. Under degeneracy the
working set can cycle, and Disser and Mosis recently showed that the exponential path survives the
generalization in a strong form: there is a polynomial objective of degree $n$ on which the
active-set method, started at the origin, visits all $2^n$ vertices of the hypercube whatever the
pivot rule \cite{dissermosis2025}.

\subsection{Interior-point methods}\label{ssec:h-ipm}
Khachiyan's ellipsoid method placed linear programming in polynomial time \cite{khachiyan1979}, and
Karmarkar's projective method began the interior-point era \cite{karmarkar1984}, whose
path-following methods need $O(\sqrt N\,\log(1/\eps))$ iterations for $N$ inequality constraints
\cite{wright1997,renegar2001}. These methods acquired their own constructions. Todd and Ye showed
that long-step primal--dual methods may need at least $n^{1/3}$ iterations to reduce the duality gap
by a constant factor \cite{toddye1996}. Deza, Nematollahi, Peyghami, and Terlaky showed that
redundant constraints can force the central path through small neighbourhoods of all $2^n$ vertices
of a Klee--Minty cube \cite{deza2006}, so that path-following methods need a number of iterations
exponential in the dimension~$n$, although not in the number of constraints
\cite{deza2008,deza2009curvature}. Allamigeon, Benchimol, Gaubert, and Joswig used tropical geometry
to show that log-barrier interior-point methods are not strongly polynomial \cite{allamigeon2018}, a
result later extended to every self-concordant barrier \cite{allamigeon2022} and to short-step
methods \cite{zong2024}.

\subsection{Descent methods: gradient and Newton}\label{ssec:h-descent}
A second tradition places the worst case against a \emph{class} of methods rather than a single rule.
Nemirovski and Yudin's information-based complexity bounds the number of oracle calls that any
method in a class needs on a class of functions, by building the hard function adaptively against
the method \cite{nemirovskiyudin1983}. Nesterov's ``worst function in the world'', an ill-conditioned
quadratic with a chain structure, gives the matching lower bounds for first-order methods on smooth
convex problems and shows that his accelerated method is optimal \cite{nesterov1983,nesterov2004}.
The programme continues for nonconvex problems \cite{carmon2020} and for stochastic first-order
oracles \cite{arjevani2023}, and the performance estimation problem computes exact worst cases of
first-order methods numerically, by semidefinite programming \cite{drori2014,taylor2017}. These
lower bounds hold for every method of the class; the constructions surveyed here are instead built
against a specific rule, and a different rule of the same class typically defeats them. The
ill-conditioned quadratic belongs to both traditions.

Newton's method has worst cases of both kinds as well. Near a minimizer with nonsingular Hessian it
converges quadratically, but at a degenerate minimizer, where the Hessian is singular, the convergence is only linear
\cite{deckerkelley1980}; Griewank's survey in this journal reviews the theory and the remedies
\cite{griewank1985}. Far from a solution the picture is
no better: Cartis, Gould, and Toint constructed functions on which Newton's method needs almost
$\eps^{-2}$ iterations to reduce the gradient norm below $\eps$, as many as steepest descent
\cite{cartis2010}. Cubic regularization improves the guarantee to $\eps^{-3/2}$
\cite{nesterovpolyak2006}, which Cartis, Gould, and Toint also show to be tight, and oracle lower
bounds for second-order methods are known in the convex case \cite{arjevani2019}.

\subsection{Tree search and local search: branch-and-bound}\label{ssec:h-tree}
Branch-and-bound \cite{landdoig1960} has its own constructions: integer programs whose linear
relaxations give bounds too weak to prune, so that exponentially many nodes must be explored
\cite{jeroslow1974,chvatal1980}. Recent work has made such lower bounds robust: they hold for
branch-and-bound trees that branch on general disjunctions, not only on variables
\cite{dey2023lower}, and they have been related to the lengths of cutting-plane proofs
\cite{basu2023}. Local search was given a complexity theory of its own with the class PLS of
polynomial local search problems \cite{johnson1988}, together with natural problems on which local
search can take exponentially many improving steps \cite{schafferyannakakis1991}. Concrete
exponential examples followed for widely used heuristics: Lloyd's $k$-means algorithm
\cite{arthurvassilvitskii2006,vattani2011} and the 2-opt heuristic for the traveling salesman problem
\cite{englert2014}.

\subsection{Randomized search: metaheuristics}\label{ssec:h-random}
For evolutionary algorithms and other randomized search heuristics, worst cases take the form of
fitness landscapes on which the expected optimization time is exponential. Deceptive and trap
functions, whose local slope points away from the global optimum, were studied early
\cite{deb1993}, and a rigorous runtime analysis developed from the late 1990s onward, including
functions on which even the simplest evolutionary algorithm needs exponential expected time
\cite{droste2002} and a set of drift-analysis tools for such bounds
\cite{doerr2012drift,lehrewitt2014}; see \cite{auger2011,jansen2013,doerrneumann2020} for
monographs. The no-free-lunch theorems frame the limits of black-box search in general
\cite{wolpertmacready1997}, and the overlap-gap property explains the failure of local algorithms on
certain random instances \cite{gamarnik2021ogp}.

\subsection{Beyond the worst case, and what the strands share}\label{ssec:h-beyond}
In practice the worst cases of the simplex method are almost never met, and a further line of work
explains why its worst-case bounds are so pessimistic. Borgwardt and Smale showed that variants of the simplex
method need polynomially many pivots on average under natural probability models
\cite{borgwardt1987,smale1983}. Spielman and Teng's smoothed analysis went further: the
shadow-vertex simplex method needs polynomially many pivots in expectation once the data of an
\emph{arbitrary} instance are perturbed by Gaussian noise, with a bound that is polynomial in the
problem size and in the reciprocal of the noise level \cite{spielmanteng2004}. The bounds have been
sharpened repeatedly \cite{dadushhuiberts2020,huiberts2025}, most recently to an optimal dependence
on the noise level \cite{bachhuiberts2025}, and smoothed analysis has become one of the main
frameworks for analysing algorithms beyond the worst case \cite{roughgarden2021}. Not every
construction yields to it: an exponential lower bound for branch-and-bound persists under Gaussian
noise \cite{dey2023lower}. We return to the robustness of the constructions in \cref{sec:boundaries}.

Each strand built its constructions with its own tools and vocabulary: polytope combinatorics for
pivoting, barrier geometry for interior-point methods, resisting oracles for first-order methods,
integrality gaps for branch-and-bound, fitness landscapes for search heuristics. The rest of this
survey collects the rule-specific constructions across these strands and asks what they have in
common. The answer we propose is a shared form: in every case the
optimum lies close to the start, and the method's own rule is misled onto a long detour or a loop.

\section{A common framework: rules on state spaces}\label{sec:framework}
The seven worst cases live in visibly different worlds---vertices of a polytope, iterates in
$\R^n$, nodes of a search tree, populations of candidate solutions. To treat them as one family we
strip each algorithm to the same three parts: the \emph{states} it moves among, the \emph{rule} by
which it moves, and a \emph{potential} that measures how much work remains.

\begin{definition}\label{def:framework}\dname{Algorithm as a rule on a state space}
An \emph{instance} of an algorithm $A$ is a tuple $I=(S,s_0,\Sstar,T,\Phi)$ where $S$ is the
\emph{state space}, $s_0\in S$ the start, $\Sstar\subseteq S$ the \emph{target} states (an optimum
or a certificate of one), $T$ the \emph{transition rule} carrying a state to its successor, and
$\Phi:S\to\R_{\ge0}$ a \emph{progress potential} with $\Phi(s)=0$ iff $s\in\Sstar$. At each state
$s$ a nonempty set $\mathcal A(s)\subseteq S$ of \emph{admissible} successors encodes the moves the
algorithm class permits; the rule selects $T(s)\in\mathcal A(s)$.
\end{definition}

The potential is the ``distance to done'': the distance to the minimizer for a descent method, the
graph distance to the optimal vertex on the polytope skeleton for pivoting, the amount of tree
still to be explored before optimality is certified for branch-and-bound, the expected number of
population members outside the basin of the optimum for a population method. The admissible set
$\mathcal A(s)$ keeps the comparison fair: it contains what a well-guided member of the \emph{same}
class could do from $s$---a different step length along the same direction, a different improving
pivot, a different branching decision---not what an outside oracle could do.

\begin{definition}\label{def:costs}\dname{Guided cost, rule cost, blow-up}
For accuracy $\eps>0$, the \emph{guided cost} $\Ng(I;\eps)$ is the least number of admissible steps
in which some admissible route drives $\Phi$ to at most $\eps$ from $s_0$; the \emph{rule cost}
$\Nr(I;\eps)$ is the number of steps the rule $T$ itself takes to do so. Their ratio is the
\emph{blow-up} $\Bl(I;\eps)=\Nr(I;\eps)/\Ng(I;\eps)$.
\end{definition}

The guided cost $\Ng$ is the length of the short route on the instance itself. It need not be a
constant: it may grow polynomially with the dimension or with $\log(1/\eps)$. A \Uex\ requires only
that the rule cost $\Nr$ grows much faster than $\Ng$. Both costs refer to the same instance. A
third quantity is used for comparison only: the \emph{baseline cost} $\Nz$, the cost of the same
rule on the \emph{normal} instance, that is, the instance with the control parameter switched off
(a quadratic bowl instead of a flat one, an undeformed cube, a single basin).

\begin{definition}\label{def:misdirection}\dname{Misdirection}
At a state $s_k$ visited by the rule, let
\[
  \Delta^{\mathrm g}_k=\Phi(s_k)-\min_{s'\in\mathcal A(s_k)}\Phi(s')
  \qquad\text{and}\qquad
  \Delta^{\mathrm r}_k=\Phi(s_k)-\Phi(T(s_k))
\]
be the largest progress any admissible step can make and the progress the rule actually makes. The
\emph{misdirection} of the step is
\[
  \delta_k=1-\frac{\max\{0,\Delta^{\mathrm r}_k\}}{\Delta^{\mathrm g}_k}\in[0,1].
\]
\end{definition}

Thus $\delta_k=0$ is the best case: the rule moves to the admissible successor with the least
remaining potential, so no step could have done better. At the other extreme, $\delta_k=1$ means
that the step makes no progress or moves backward. Intermediate values measure the share of the
available progress that the step wastes. For example, at the start of the Klee--Minty cube of
\cref{ssec:simplex} the optimum is an adjacent vertex, so $\Delta^{\mathrm g}_0=1$; Dantzig's rule
moves instead to a vertex at distance two, so $\Delta^{\mathrm r}_0=-1$ and $\delta_0=1$. On a
normal instance the rule is well guided and $\delta_k$ is zero or small; a \Uex\ is an instance on
which it is large for long enough to matter.

Misdirection is measured through the potential $\Phi$, not through an angle between the step and
the direction of the optimum, so it applies verbatim to pivoting, tree search, and populations,
where such an angle has no meaning. \Cref{tab:framework} instantiates the framework for each member.

\begin{table}[t]
\renewcommand{\arraystretch}{1.25}\footnotesize\centering
\caption{One framework, seven instantiations. Every member is a rule on a state space with a
progress potential and a set of admissible moves. The last column gives the misdirection $\delta_k$
(\cref{def:misdirection}) \emph{on the member's \Uex}, where the rule is misdirected by
construction; on the corresponding normal instance $\delta_k$ is zero or small.}
\label{tab:framework}
{\setlength{\tabcolsep}{4.5pt}\begin{tabular}{@{}p{1.85cm} >{\raggedright\arraybackslash}p{1.7cm} >{\raggedright\arraybackslash}p{2.1cm} >{\raggedright\arraybackslash}p{2.35cm} >{\raggedright\arraybackslash}p{3.45cm}@{}}
\toprule
\textbf{Member} & \textbf{State} & \textbf{Potential $\Phi$} & \textbf{Admissible moves} & \textbf{Misdirection on the \Uex} \\
\midrule
Newton & iterate $x\in\R^n$ & distance $\|x-x^\star\|$ & steps of any length along the Newton direction & $\delta_k=\tfrac{2m-2}{2m-1}$ at every step: the step is too short by a factor $2m-1$ \\
Gradient & iterate $x\in\R^n$ & distance $\|x-x^\star\|_\infty$ & steps of any length along $-\nabla f$ & $\delta_k\approx\tfrac{\rho}{1+\rho}$ at every step: the step is too short by a factor of about $1+\rho$ \\
Simplex & vertex of the polytope & graph distance to $x^\star$ & improving pivots & $\delta_k=1$ on $2^{n-1}-1$ pivots, which move away from $x^\star$; $\delta_k=0$ on the others \\
Interior-point & interior point & duality gap & interior steps in the feasible region & every step is confined to a neighbourhood of the bent central path \\
Branch-and-bound & set of open nodes & nodes still needed to close the tree & branching on any integer disjunction & $\delta_k=1$: each variable branch leaves two subproblems as hard as their parent \\
Metaheuristic & population in $\R^n$ & expected number of members outside the optimal basin & any placement of a member & defined in expectation: a member lands in the optimal basin with probability $k^{-n}$ \\
Active-set & working set at a point & distance $\|x-x^\star\|$ & any change of the working set & $\delta_k=1$: zero-length steps change the working set but not the point \\
\bottomrule
\end{tabular}}
\end{table}

\begin{definition}\label{def:uex}\dname{\Uex}
A family $\{I_N\}_{N\ge1}$ of instances for a rule $T$ is a \emph{\Uex} if
\emph{(i)} the guided cost stays small, $\Ng(I_N;\eps)\le p(N,\log1/\eps)$ for a fixed polynomial
$p$, and \emph{(ii)} the blow-up is unbounded, $\Bl(I_N;\eps)\to\infty$ as $N\to\infty$ for
each fixed $\eps\in(0,\Phi(s_0))$, including the degenerate case $\Nr=\infty$ in which the rule never
reaches the target. The index $N$ is the instance's \emph{control parameter}---informally, the knob
that turns the worst case on: a dimension, a flatness order, a number of basins. For a randomized
rule, $\Nr$ is replaced by its expectation.
\end{definition}

Misdirection and blow-up are linked by a simple accounting identity. Let the rule visit
$s_0,s_1,\dots,s_{\Nr}$ with $\Phi(s_{\Nr})\le\eps$, and assume that every non-target state has an
admissible successor with strictly smaller potential, so that $\Delta^{\mathrm g}_k>0$ and
$\delta_k$ is defined. Let $\Delta^\ast$ be an upper bound on the available progress
$\Delta^{\mathrm g}$ at every state visited by the rule or by a shortest guided route, and write
\[
  R=\sum_{k<\Nr}\max\{0,\Phi(s_{k+1})-\Phi(s_k)\},\qquad
  \bar\delta=\frac1{\Nr}\sum_{k<\Nr}\delta_k,\qquad
  \kappa=\frac{\Ng(I;\eps)\,\Delta^\ast}{\Phi(s_0)-\eps}
\]
for the rule's total \emph{regress} (potential it gives back), its mean misdirection, and the
\emph{efficiency} of the guided route ($\kappa\ge1$, with equality when every guided step removes the
largest possible progress).

\begin{proposition}\label{prop:blowup}\dname{Two sources of blow-up}
Under the assumptions above,
\[
  \Nr(I;\eps)\ \ge\ \frac{\Phi(s_0)-\eps+R}{\Delta^\ast\,(1-\bar\delta)},
  \qquad\text{hence}\qquad
  \Bl(I;\eps)\ \ge\ \frac{1+R/(\Phi(s_0)-\eps)}{\kappa\,(1-\bar\delta)}.
\]
\end{proposition}

\begin{figure}[t]\centering
\includegraphics[width=\textwidth]{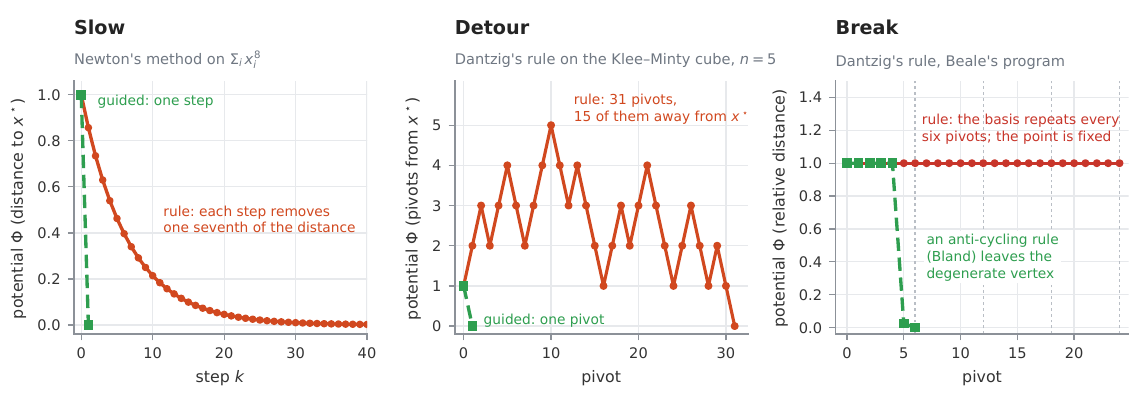}
\caption{The three ways a rule is forced the long way, shown as the potential $\Phi$ along computed
trajectories. \emph{Slow:} Newton's method on the flat-bottomed bowl $\sum_i x_i^{8}$ removes one
seventh of the remaining distance per step; the guided route (green) takes one step.
\emph{Detour:} Dantzig's rule on the Klee--Minty cube with $n=5$ starts one pivot from the optimum
and takes $31$, of which $15$ move away from it. \emph{Break:} on Beale's degenerate linear program
\cite{beale1955}, Dantzig's rule returns to its starting basis every six pivots and the point never
moves; Bland's anti-cycling rule (green) leaves the degenerate vertex and reaches the optimum in six
pivots.}
\label{fig:mechanisms}
\end{figure}

The proof sums the progress of the individual steps and is given in \cref{app:proofs}. The bound names the only two ways a rule can be forced the long way, and \cref{fig:mechanisms} shows each on
a computed trajectory.\footnote{Code that reproduces Figures~\ref{fig:mechanisms}, \ref{fig:remedy}, and~\ref{fig:fragility} is provided as supplementary material.} Either the rule captures only a small share of the progress on offer
($\bar\delta\to1$): it is \emph{slow}, as Newton's method is on a flat bottom. Or it repeatedly
gives back potential it has already gained ($R\gg\Phi(s_0)$): it \emph{detours}, as the simplex
method does on the Klee--Minty cube, where the bound is attained (\cref{ssec:simplex}). A
rule that returns to a state it has already visited never terminates: around the cycle its gains are
cancelled by regress, or---as in degenerate pivoting, where $\Phi$ does not move at all---every step
has $\delta_k=1$. Either way $\Nr=\infty$, and the rule \emph{breaks}. Turning these mechanisms into
exact growth in each member's control parameter and in the dimension $n$ is the work of
\cref{sec:members,sec:theorem}.

\section{The seven members}\label{sec:members}
\begin{figure}[p]\centering
\includegraphics[width=0.9\textwidth]{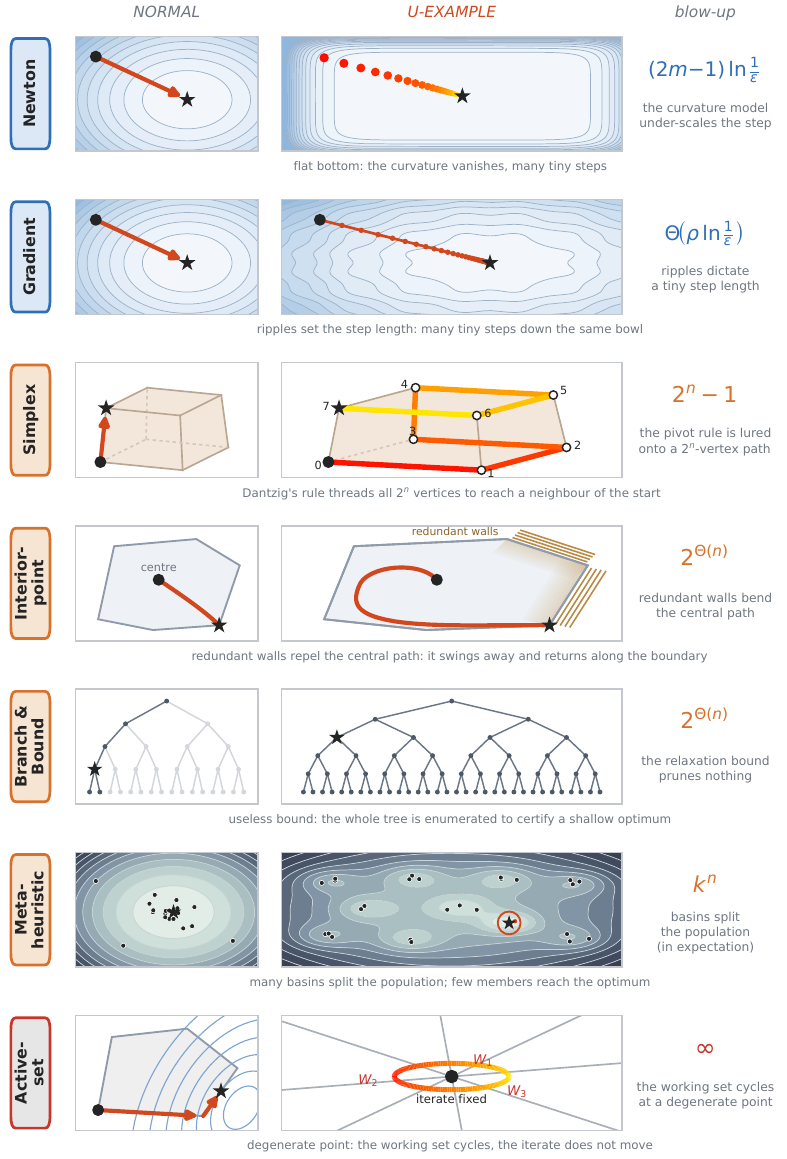}
\caption{The family in one plate. Each row is one algorithm class, in
descent$\to$search$\to$break order (matching \cref{fig:story}), with its dimension class on the
left spine (blue: dimension-free; orange: exponential; red: break). For every member a \emph{normal}
instance (left) sits beside its \Uex\ (right), drawn with the algorithm's trajectory---%
Newton's method crawling into a flat bottom, gradient descent stuttering on a rippled bowl, Dantzig's
rule threading all $2^n$ vertices of the Klee--Minty cube ($0\to7$) to a vertex adjacent to the
start, the interior-point path bent by redundant walls, branch-and-bound enumerating a tree it
cannot prune, a population fragmenting across many basins, and an active-set method changing its
working set at a degenerate point without moving. The right column gives each blow-up
\cite{griewank1985,kleeminty1972,deza2008,jeroslow1974,deb1993,nocedalwright2006}. The planar
panels are schematic; the costs are derived in the text.}
\label{fig:plate}
\end{figure}

We treat seven representative classes; \cref{ssec:beyond} points to the worst cases of other
methods. For each member we review the background, state the \Uex, compute its cost in general
dimension~$n$, and summarize the state of the art. \Cref{fig:plate} draws all seven in one visual
grammar: a normal instance beside its \Uex, with the algorithm's trajectory on each.

\paragraph{Baseline principle} Three costs appear for every member (\cref{sec:framework}): the
baseline cost $\Nz$ of the rule on the normal instance, and, on the \Uex\ itself, the guided cost
$\Ng$ and the rule cost $\Nr$. The blow-up is the ratio $\Nr/\Ng$ of the last two. We state each
cost as it is, including a guided cost that grows with $n$ or with $\log(1/\eps)$, rather than
idealize it to a constant.

\subsection{Newton's method (flat bottom)}\label{ssec:newton}
\emph{Background.} Newton's method minimizes a local quadratic model at each step and converges
quadratically near a minimizer at which the Hessian is nonsingular
\cite{nocedalwright2006,boydvandenberghe2004}. When the Hessian is singular at the minimizer the
quadratic model is a poor guide, and the convergence is only linear. In one variable the picture is
classical: at a root of multiplicity $p$ of the equation being solved, the error is multiplied by
$1-1/p$ at each step, and lengthening the step by the factor $p$ restores fast convergence. The
theory in several variables begins with the analysis of Newton's method at singular points
\cite{deckerkelley1980} and is reviewed, together with the remedies, in Griewank's survey
\cite{griewank1985}. Global worst cases are known too: there are functions on which
Newton's method needs almost $\eps^{-2}$ iterations to drive the gradient norm below $\eps$
\cite{cartis2010}.

\emph{The \Uex.} Take the separable flat-bottomed bowl $f(x)=\sum_{i=1}^n x_i^{2m}$ with $m\ge1$,
whose minimizer is $x^\star=0$. The Newton step is the same in every coordinate,
\[
  x_{k+1}=x_k-\nabla^2f(x_k)^{-1}\nabla f(x_k)=\frac{2m-2}{2m-1}\,x_k .
\]
For $m=1$ the bowl is quadratic and the method lands on the minimizer in one step
(\cref{fig:plate}, Newton row, left panel). For $m\ge2$ the bottom is flat, the model sees almost no
curvature, and the iterate removes only the fraction $1/(2m-1)$ of the remaining distance at each
step (right panel). The short route is one step long: the same Newton direction, with step length
$2m-1$ instead of $1$, reaches $x^\star$ exactly.

\emph{Cost.} Take the distance $\Phi(x)=\|x-x^\star\|$ as potential and steps of any length along
the Newton direction as admissible moves. \textbf{Control parameter:} the flatness order $2m$.
\textbf{Baseline:} $\Nz=1$ on the quadratic. \textbf{Guided cost:} $\Ng=1$. \textbf{Rule cost:}
for $m\ge2$, since each step multiplies the distance by $(2m-2)/(2m-1)$,
\[
  \Nr=\Bigl\lceil \ln\tfrac{\Phi(x_0)}{\eps}\Big/\ln\tfrac{2m-1}{2m-2}\Bigr\rceil,
  \qquad (2m-2)\ln\tfrac{\Phi(x_0)}{\eps}\ \le\ \Nr\ \le\ (2m-1)\ln\tfrac{\Phi(x_0)}{\eps}+1 .
\]
\begin{payload}\centering
$\displaystyle \Bl_{\mathrm{Newton}}(m,n;\eps)\approx(2m-1)\ln\frac{\Phi(x_0)}{\eps}$\\[3pt]
(linear in the flatness order; independent of $n$).\par
\end{payload}
Every step has the same misdirection $\delta_k=(2m-2)/(2m-1)$, and there is no regress: this is
the \emph{slow} mechanism of \cref{prop:blowup} in pure form. Because the coordinates decouple, the
dimension leaves the count unchanged.

\begin{remark}\label{rem:newton-f}
The choice of potential matters here. Measured in function value, the slowdown disappears:
$f(x_k)=f(x_0)\,r^{2mk}$ with $r=(2m-2)/(2m-1)$, and $r^{2m}\to e^{-1}$ as $m$ grows, so about
$\ln(f(x_0)/\eps)$ steps reduce $f$ below $\eps$ for every $m$. A flat bottom makes the function
value a poor measure of progress: points far from the minimizer already have nearly optimal values.
The \Uex\ concerns accuracy in the solution, not in the objective.
\end{remark}

\emph{State of the art.} When the multiplicity is known, the lengthened step restores fast
convergence; \cite{griewank1985} reviews this and other modifications. When it is not, regularized Newton methods give global guarantees
without it: cubic regularization \cite{nesterovpolyak2006}, and, more recently, regularization
proportional to a power of the gradient norm, which keeps the cost of a Newton step
\cite{mishchenko2023,doikov2024}. Oracle lower bounds for second-order methods bound what any such
scheme can achieve \cite{arjevani2019}.

\subsection{Gradient descent (rippled bowl)}\label{ssec:gradient}
\emph{Background.} Gradient descent with the standard step length $1/L$, where $L$ is the Lipschitz
constant of the gradient, needs $\Theta(\kappa\log(1/\eps))$ iterations on a quadratic with
condition number $\kappa$. This is the classical worst case of the method: the iterate zig-zags
across a narrow valley. Nesterov's worst-case function, an ill-conditioned quadratic built on the
discrete Laplacian, turns it into a lower bound of order $\sqrt\kappa\,\log(1/\eps)$ for \emph{all}
first-order methods, attained by acceleration
\cite{nesterov1983,nesterov2004,nesterov2018,nemirovskiyudin1983,bubeck2015}. In that example the
difficulty is intrinsic to first-order information. The \Uex\ below isolates a different and purely
rule-specific difficulty: the step length $1/L$ is dictated by the largest curvature anywhere, even
where that curvature is irrelevant to the route.

\emph{The \Uex.} Add fine ripples to a perfectly conditioned bowl:
\begin{equation}\label{eq:rippled}
  f(x)=\sum_{i=1}^n\Bigl(\tfrac12x_i^2+\frac{\rho}{\omega^2}\cos(\omega x_i)\Bigr),
  \qquad \rho\ge1,\ \ \omega>0 .
\end{equation}
The ripples have height $\rho/\omega^2$ and slope at most $\rho/\omega$, both negligible when the
frequency $\omega$ is large, but their curvature is $\rho$ whatever the frequency: the second
derivative in each coordinate is $1-\rho\cos(\omega x_i)\in[1-\rho,\,1+\rho]$, so that $L=1+\rho$.
Seen from any distance the function is the bowl $\tfrac12\|x\|^2$, on which gradient descent with
unit step reaches the minimizer in one step (\cref{fig:plate}, gradient row, left panel). The
ripples force the step length down to $1/(1+\rho)$, and the iterate descends the same bowl in many
small steps (right panel).

\emph{Cost.} Take $\Phi(x)=\|x\|_\infty$, the start $x_0=(1,\dots,1)$, and a target accuracy $\eps$
with $2\rho/\omega\le\eps<1$, so that the ripples are finer than the accuracy sought. Admissible
moves are steps of any length along $-\nabla f$. \textbf{Control parameter:} the roughness $\rho$.
\textbf{Baseline:} $\Nz=1$ on the smooth bowl ($\rho=0$). \textbf{Guided cost:} $\Ng=1$, because
the unit step from $x_0$ lands within $\rho/\omega\le\eps/2$ of the minimizer. \textbf{Rule cost:}
by \cref{lem:rippled} in \cref{app:proofs},
\[
  \tfrac{2\rho-1}{3}\,\ln\tfrac1\eps\ \le\ \Nr\ \le\ 2(1+\rho)\ln\tfrac1\eps+1 .
\]
\begin{payload}\centering
$\displaystyle \Bl_{\mathrm{gradient}}(\rho,n;\eps)=\Theta\bigl(\rho\,\ln\tfrac1\eps\bigr)=\Theta\bigl(L\ln\tfrac1\eps\bigr)$\\[3pt]
(linear in the roughness; independent of $n$).\par
\end{payload}
In a direct computation the count is close to $(1+\rho)\ln(1/\eps)$: with $\eps=10^{-3}$ the method
takes $114$ steps for $\rho=16$ and $1771$ for $\rho=256$. Every step has misdirection about
$\rho/(1+\rho)$: as for Newton's method, the direction is right and the length is wrong.

The rippled bowl differs from the ill-conditioned quadratic in kind. The quadratic is slow for
every first-order method, and acceleration improves the count only from $\kappa$ to $\sqrt\kappa$.
The rippled bowl is slow only for the fixed step $1/L$: any rule that tests a longer step, such as
a line search, solves it in one step. We have not found this instance in the literature, but its
mechanism is familiar, and the two examples are complementary: for natural families of
ill-conditioned quadratics, such as discretized Laplacians, the condition number grows
polynomially with the dimension, whereas the roughness of \cref{eq:rippled} is independent of it.

\begin{remark}
The instance is a \Uex\ only down to the scale of the ripples. Below the accuracy $\rho/\omega$ the
ripples create spurious stationary points, at which gradient descent stops: a break, not a slow
crawl. This is the mechanism of the metaheuristic member (\cref{ssec:meta}), in one dimension.
\end{remark}

\emph{State of the art.} Step-length rules that adapt to the local rather than the global curvature
are the classical answer, from line searches \cite{nocedalwright2006} to the two-point step sizes
of Barzilai and Borwein \cite{barzilaiborwein1988}. The same tension between step length and
sharpest curvature has resurfaced in machine learning, where full-batch gradient descent on neural
networks is observed to operate at the \emph{edge of stability}, with the largest Hessian eigenvalue
hovering near $2$ divided by the step length \cite{cohen2021}. Exact worst-case rates of first-order
methods are now computed routinely by performance estimation \cite{drori2014,taylor2017}.

\subsection{Simplex method (tilted cube)}\label{ssec:simplex}
\emph{Background.} The simplex method walks along the edges of the feasible polytope, choosing at
each vertex an improving edge by a pivot rule \cite{dantzig1963}. Its efficiency in practice has
always contrasted with its worst case. Every classical exponential example, from Klee and Minty's
onward \cite{kleeminty1972,jeroslow1973,avischvatal1978,goldfarbsit1979}, uses the same template, a
\emph{deformed product} \cite{amentaziegler1999}: a polytope combinatorially equivalent to a cube,
with facets tilted so that a monotone edge path through all $2^n$ vertices exists, and scaled so
that a given rule follows it. Each construction is tied to one rule, and newer rules have required
new constructions (\cref{sec:history}).

\emph{The \Uex.} We use the Klee--Minty linear program in its classical scaled form,
\begin{equation}\label{eq:km}
  \max\ \sum_{j=1}^n 2^{\,n-j}x_j\quad\text{subject to}\quad
  2\sum_{j<i}2^{\,i-j}x_j+x_i\le5^{\,i}\ \ (i=1,\dots,n),\qquad x\ge0 ,
\end{equation}
with Dantzig's largest-coefficient rule started at the vertex $x=0$. The feasible region is
combinatorially an $n$-cube: each vertex is labelled by $b\in\{0,1\}^n$, with $b_i=1$ when $x_i$ is
basic, and two vertices are adjacent exactly when their labels differ in one coordinate. The optimum
is $x^\star=(0,\dots,0,5^n)$, with label $(0,\dots,0,1)$: it is \emph{adjacent} to the start, and
the single pivot that brings $x_n$ into the basis solves the problem (the normal panels of the
simplex rows in \cref{fig:plate,fig:plate3d}). Dantzig's rule instead enters $x_1$, whose objective
coefficient is largest, and goes on to visit all $2^n$ vertices before it arrives
\cite{kleeminty1972}. We confirmed this by direct computation for $n\le12$: the method takes
exactly $2^n-1$ pivots, and the labels of the vertices it visits follow the reflected Gray code. That
such a path can be monotone at all is the geometric content of the construction;
\cref{lem:km} in \cref{app:proofs} gives a short proof for the deformed cube
$\{\,0\le x_1\le1,\ \tau x_{j-1}\le x_j\le1-\tau x_{j-1}\,\}$ with objective $x_n$, the form
studied in \cite{gartnerhenkziegler1998}.

\emph{Cost.} \textbf{Potential:} the graph distance to $x^\star$ on the skeleton of the cube, which
is the Hamming distance between labels; admissible moves are improving pivots. Every non-optimal
vertex has an improving neighbour one step closer to the optimum, so the available progress is
$\Delta^{\mathrm g}=1$ at every state. \textbf{Control parameter:} the dimension $n$, together
with the scaling of \cref{eq:km}. \textbf{Baseline:} $\Nz=1$ on the undeformed cube.
\textbf{Guided cost:} $\Ng=1$. \textbf{Rule cost:} $\Nr=2^n-1$. Each pivot changes the distance by
exactly $\pm1$, and the walk starts at distance $1$ and ends at $0$; hence exactly $2^{n-1}$ pivots
move toward the optimum and the remaining $R=2^{n-1}-1$ move away from it.
\begin{payload}\centering
$\displaystyle \Bl_{\mathrm{simplex}}(n)=2^n-1$\\[3pt]
(exponential in $n$).\par
\end{payload}
The misdirection here is regress, not slowness: $\delta_k$ is $1$ on the pivots that move away
and $0$ on the others, so the mean is $\bar\delta=(2^{n-1}-1)/(2^n-1)\approx\tfrac12$, and the rule
wastes about half its pivots walking away from an optimum that was one pivot from the start. With
$\Delta^\ast=1$, the bound of \cref{prop:blowup} is attained in the limit $\eps\to0$:
$\Nr(1-\bar\delta)=2^{n-1}=\Phi(s_0)+R$.

\emph{State of the art.} The instance is tuned to one rule and one scaling. In the same
computation, the steepest-edge and greatest-improvement rules solve \cref{eq:km} in a single pivot;
a uniformly random improving pivot needs about $20$ pivots on average at $n=12$; and Bland's rule
needs $465$. Even Dantzig's rule takes only $2n-1$ pivots once each constraint of \cref{eq:km} is
rescaled to unit norm, which leaves the linear program unchanged. Each of these rules has its own
exponential instance, however: greatest improvement \cite{jeroslow1973}, Bland \cite{avischvatal1978},
steepest edge \cite{goldfarbsit1979}, the history-based rules of Zadeh and Cunningham
\cite{disser2023zadeh,avisfriedmann2017}, and, most recently, every norm of steepest edge and every
projection of the shadow-vertex rule \cite{black2026}. Random-facet rules admit subexponential upper
and lower bounds \cite{kalai1992,matousek1996,friedmannhansenzwick2011}. Whether some pivot rule is
polynomial remains open, as does the polynomial Hirsch conjecture, which any such rule would imply
\cite{santos2012,kalaikleitman1992,todd2014}. In the other direction, smoothed analysis shows that
the expected number of pivots of the shadow-vertex rule is polynomial once the data are randomly
perturbed \cite{spielmanteng2004,dadushhuiberts2020,huiberts2025,bachhuiberts2025}; we return to
this in \cref{sec:boundaries}.

\subsection{Interior-point methods (redundant walls)}\label{ssec:ipm}
\emph{Background.} Path-following interior-point methods trace the \emph{central path}, the curve
of minimizers of the objective plus a decreasing multiple of a logarithmic barrier, from the
analytic centre of the feasible region to an optimal face
\cite{karmarkar1984,wright1997,renegar2001}. Their iteration bound, $O(\sqrt N\log(1/\eps))$ for
$N$ inequality constraints, depends on the description of the feasible region and not only on its
shape: the barrier has one term per constraint, redundant or not. Lower bounds exploit this. The
central path can be made to bend sharply \cite{deza2009curvature}, and the number of sharp turns
bounds the number of iterations of any method that follows it.

\emph{The \Uex.} Deza, Nematollahi, and Terlaky \cite{deza2008}, refining \cite{deza2006}, start
from the Klee--Minty cube in dimension $n$ and add, for each of the $n$ facets through the origin,
many copies of a constraint parallel to that facet and slightly outside the cube. The copies are
redundant: the feasible region is still the cube, described by its original $2n$ inequalities. But
each copy adds a term to the barrier, and a barrier term repels: the more copies stand behind a
facet, the harder the central path is pushed away from it. Together the copies drive the path
against the rest of the boundary, so that it passes within a small distance $\delta$ of every one
of the $2^n$ vertices before it reaches the optimum (\cref{fig:plate}, interior-point row, where
the path swings away from the walls and returns along the boundary). The construction uses
$N=O(2^{2n}n^3)$ inequalities in all.

\emph{Cost.} Take the duality gap as potential and interior steps in the feasible region as
admissible moves. \textbf{Control parameter:} the redundancy, that is, the numbers of parallel
copies, with the dimension $n$. \textbf{Baseline:} on the cube without redundant constraints,
$\Nz=O(\sqrt n\,\log(1/\eps))$. \textbf{Guided cost:} the redundant constraints do not change the
feasible region, so the central path of the irredundant description is an admissible interior
route, and $\Ng=O(\sqrt n\,\log(1/\eps))$. \textbf{Rule cost:} a path-following method must turn
with the central path, which gives $\Nr\ge2^n-1$ \cite{deza2008}.
\begin{payload}\centering
$\displaystyle \Bl_{\mathrm{IPM}}(n;\eps)\ \ge\ \frac{2^n-1}{O(\sqrt n\,\log(1/\eps))}=2^{\Theta(n)}$\\[3pt]
(exponential in the dimension $n$, for fixed $\eps$).\par
\end{payload}
The statement needs care, and the framework makes the care easy to express. The blow-up is
exponential in the \emph{dimension}, not in the size of the input. Measured against the number of
constraints, $2^n$ iterations is of order $\sqrt N$ up to polynomial factors in $n$, which is what
the upper bound $O(\sqrt N\log N)$ allows: the instance shows that the classical iteration bound is
nearly tight, and interior-point methods remain polynomial. The misdirection is real nonetheless.
The optimum and a short interior route to it are the same with and without the redundant walls;
only the rule's reading of the region has changed. A planar region admits only a modest bend, so
the interior-point row of \cref{fig:plate} understates the effect; the algebra is that of
\cite{deza2008}.

\emph{State of the art.} Redundant constraints are removed by presolving in practice. In theory,
lower bounds of the same flavour have been pushed much further: log-barrier methods are not
strongly polynomial \cite{allamigeon2018}, nor is any method based on a self-concordant barrier
\cite{allamigeon2022} or any short-step method \cite{zong2024}; and long-step methods need at
least order $n^{1/3}$ iterations to reduce the duality gap by a constant factor \cite{toddye1996}.

\subsection{Branch-and-bound (unprunable tree)}\label{ssec:bb}
\emph{Background.} Branch-and-bound solves an integer program by recursive splitting
\cite{landdoig1960,wolsey1998,conforti2014ip}. At each node it solves the linear relaxation; it
discards the node if the relaxation is infeasible, has an integral solution, or has a bound no
better than the best solution known, and otherwise it branches, usually on a variable with a
fractional value. The trajectory is the sequence of nodes processed, and the method is only as
good as its bounds: if the relaxation never rules anything out, the tree is enumerated.

\emph{The \Uex.} Jeroslow's integer program \cite{jeroslow1974} is the cleanest instance:
\begin{equation}\label{eq:jeroslow}
  \max\ x_1\quad\text{subject to}\quad 2x_1+2x_2+\dots+2x_n=n,\qquad x\in\{0,1\}^n ,
\end{equation}
with $n$ odd. The left-hand side is even and the right-hand side is odd, so the program is
infeasible, and a one-line parity argument certifies it. Branch-and-bound on variables cannot see
this. After any $d\le(n-1)/2$ variables have been fixed, the remaining ones can still be set
fractionally to satisfy the equation, so the relaxation at every node of depth at most $(n-1)/2$ is
feasible and non-integral, and all of these nodes must be branched (\cref{lem:jeroslow} in
\cref{app:proofs}). The short route is a different branching decision at the root: the disjunction
$\sum_ix_i\le(n-1)/2$ or $\sum_ix_i\ge(n+1)/2$ is valid for every integer point, and both of its
sides contradict the equation, so the tree closes with three nodes.

\emph{Cost.} The state is the set of open nodes, the potential is the number of nodes still
needed to close the tree, and the admissible moves are branchings on any integer disjunction.
\textbf{Control parameter:} the parity of $n$, a switch. \textbf{Baseline:} for even $n$ the
relaxation at the root has an integral vertex, and $\Nz=1$. \textbf{Guided cost:} $\Ng=3$.
\textbf{Rule cost:} $\Nr\ge2^{(n+1)/2}-1$ nodes; a direct count gives $139$ nodes for $n=7$ and
more than $10^9$ for $n=31$.
\begin{payload}\centering
$\displaystyle \Bl_{\mathrm{B\&B}}(n)\ \ge\ \tfrac13\bigl(2^{(n+1)/2}-1\bigr)=2^{\Theta(n)}$\\[3pt]
(exponential in $n$; switched on by parity).\par
\end{payload}
Each variable branch replaces a subproblem by two that are as hard as their parent, so the work
remaining grows rather than shrinks: the misdirection is $\delta_k=1$ at every step. In
optimization versions of the phenomenon, such as Chv\'atal's knapsack problems \cite{chvatal1980},
an optimal solution is found at a shallow node and the tree is spent on certifying it, which is the
situation drawn in \cref{fig:plate}: the blow-up is in \emph{certifying} the answer, not in finding
it.

\emph{State of the art.} With arbitrary disjunctions and suitable basis reduction, integer
programming in fixed dimension is polynomial \cite{lenstra1983}, so the weakness lies in the
branching rule, not in the tree method as such. Even general disjunctions do not always help,
however: there are packing, covering, and traveling-salesman instances on which every
branch-and-bound tree with general disjunctions is exponential, and an exponential lower bound
persists when Gaussian noise is added to the data of one such family \cite{dey2023lower}. The relative power of branching and cutting
planes is now understood in some detail \cite{basu2023}.

\subsection{Metaheuristics (fragmented population)}\label{ssec:meta}
\emph{Background.} Population-based metaheuristics---evolutionary algorithms, particle swarms,
multistart local search---maintain a set of candidate solutions and move them by randomized local
rules \cite{blum2003metaheuristics}. Their worst cases are landscapes on which the local rule pulls
candidates toward the wrong place. Deceptive functions, whose slope leads away from the global
optimum, were designed for this purpose \cite{deb1993}; rigorous runtime analysis has since shown
that even simple functions can force exponential expected time on the $(1+1)$ evolutionary
algorithm, which optimizes linear functions in expected time $O(n\log n)$ \cite{droste2002}, and it
has built the drift-analysis tools that make such statements routine
\cite{doerr2012drift,lehrewitt2014,doerrneumann2020}.

\emph{The \Uex.} We state the member for the simplest model that the picture in \cref{fig:plate}
represents: a population of independent members, each of which starts at a uniformly random point
and descends to the local optimum of the basin it starts in. Let the landscape be separable on the
unit cube, $f(x)=\sum_{i=1}^n g(x_i)$, where $g$ has $k$ basins of equal width on $[0,1]$, exactly
one of which contains its global minimizer. Then $f$ has $k^n$ basins of equal volume, one of them
optimal. With $k=1$ there is a single basin and every member reaches the optimum (left panel). With
$k\ge2$ each basin captures its share of the population, and a member reaches the optimum only if it
starts in the optimal basin in every coordinate (right panel).

\emph{Cost.} The state is the population, and the potential is the expected number of members
outside the optimal basin; a guided method may place a member anywhere. \textbf{Control parameter:}
the number $k$ of basins per coordinate. \textbf{Baseline:} $\Nz=1$ descent. \textbf{Guided cost:}
$\Ng=1$. \textbf{Rule cost:} a member succeeds with probability $k^{-n}$, independently of the
others, so the expected number of members, or of restarts, until the first success is $k^n$, and a
population of size $P$ succeeds with probability $1-(1-k^{-n})^P\le P\,k^{-n}$.
\begin{payload}\centering
$\displaystyle \mathbb E\,\Bl_{\mathrm{meta}}(k,n)=k^n$\\[3pt]
(exponential in $n$; in expectation).\par
\end{payload}
Two qualifications apply. The statement is probabilistic: randomness means that no instance can
\emph{forbid} success, so the claim is that success is improbable and the expected cost huge, not
that the long way is forced. And it is a statement about the model, not about every metaheuristic:
on a separable landscape, operators that recombine coordinates of different members can assemble
the optimum from partial successes, and the general limits of black-box search are those of the
no-free-lunch theorems \cite{wolpertmacready1997}.

\emph{State of the art.} Runtime analysis gives rigorous versions of the same picture for
evolutionary algorithms on multimodal functions \cite{auger2011,jansen2013,doerrneumann2020}, and
diversity mechanisms, restarts, and niching are the standard remedies. For random landscapes, the
overlap-gap property identifies a geometric obstruction that defeats broad classes of local
algorithms \cite{gamarnik2021ogp}.

\subsection{Active-set methods (degenerate cycling)}\label{ssec:activeset}
\emph{Background.} Active-set methods solve a quadratic or nonlinear program by guessing which
inequality constraints are active at the solution. They hold a \emph{working set} of constraints as
equalities, minimize on the corresponding face, and then add a blocking constraint or drop one with
a multiplier of the wrong sign \cite{nocedalwright2006,fletcher2000}. On a strictly convex
quadratic program the method terminates finitely, provided each change of the working set is
followed by a step of nonzero length. With a linear objective the method is the simplex method, and
the working set is the basis. The class is therefore the nonlinear generalization of pivoting,
which is why its normal instance in \cref{fig:plate} has a curved objective and an optimum that
need not be a vertex.

\emph{The \Uex.} At a \emph{degenerate} point more constraints are active than are needed to
determine it, and several working sets describe the same point. The step after a change of working
set can then have zero length: the working set changes, the point does not. A sequence of such
changes can return to a working set already visited, after which the method repeats the same
sequence forever. Beale's linear program is the classical instance \cite{beale1955}: Dantzig's
rule returns to its starting basis after six pivots, without ever leaving the starting vertex
(\cref{fig:mechanisms}, right). Since the simplex method is the active-set method for a linear
objective, it is an instance for the whole class (\cref{fig:plate}, active-set row).

\emph{Cost.} The state is the working set together with the current point, the potential is the
distance to the solution, and any change of working set is admissible. \textbf{Control parameter:}
degeneracy, a switch. \textbf{Baseline:} at nondegenerate points every change of working set is
followed by a nonzero step and the method terminates finitely. \textbf{Guided cost:} finite: a
different choice of constraint to drop leaves the degenerate point. \textbf{Rule cost:}
$\Nr=\infty$.
\begin{payload}\centering
$\displaystyle \Bl_{\mathrm{active}}=\infty$\\[3pt]
(independent of the dimension: a break).\par
\end{payload}
This member differs from the simplex member in mechanism, not only in class. On the Klee--Minty
cube every pivot moves and the path is long but finite; here no step moves, and every step has
$\delta_k=1$. We keep it as the boundary of the family rather than force it into a ratio.

\emph{State of the art.} Anti-cycling rules such as Bland's \cite{bland1977}, and small
perturbations of the constraints, remove the break. The class also inherits the
\emph{detour} of the simplex method, in a strong form. Disser and Mosis construct a polynomial
objective of degree $n$ on which the active-set method started at the origin visits all $2^n$
vertices of the hypercube and takes $2^n-1$ iterations \emph{whatever the pivot rule}
\cite{dissermosis2025}---an unconditional lower bound of a kind that is not known for the simplex
method itself.

\subsection{Beyond the seven}\label{ssec:beyond}
The seven classes are representative, not exhaustive. They were chosen as the methods of the
standard curriculum, one or two for each main type of problem (\cref{sec:intro}), and
rule-specific worst cases are known for many other methods. \Cref{tab:beyond} lists a selection.
We do not work these out as members, with a potential and three costs, but each fits the
description of \cref{sec:intro}: a rule that works well in general, and an instance built against
exactly that rule. The failures are of the kinds met above. Some methods are slow, some are led
onto an exponentially long path, and a notable number break outright: they cycle, diverge, or
converge to a point that is not a solution. Several entries return in \cref{sec:connections}, and
two of them, heavy-ball momentum and cyclic coordinate descent, are examined in
\cref{sec:outlook}.

\begin{table}[t]
\renewcommand{\arraystretch}{1.25}\footnotesize\centering
\caption{Rule-specific worst cases beyond the seven members: a selection. The last column gives
the kind of failure: \emph{slow} (a poor rate), \emph{long path} (exponentially many steps), or
\emph{break} (no convergence to a solution).}
\label{tab:beyond}
\begin{tabular}{@{}>{\raggedright\arraybackslash}p{2.75cm} >{\raggedright\arraybackslash}p{6.35cm} >{\raggedright\arraybackslash}p{1.45cm} l@{}}
\toprule
\textbf{Method} & \textbf{Worst case} & \textbf{Failure} & \textbf{Source} \\
\midrule
Heavy-ball momentum & a smooth strongly convex function on which the iterates cycle & break & \cite{lessard2016,goujaud2025} \\
Cyclic coordinate descent & convex quadratics on which a fixed order is $O(n^2)$ times slower than a random one; a function of three variables on which the iterates cycle & slow; break & \cite{sunye2021,powell1973} \\
Quasi-Newton (BFGS) & nonconvex functions on which the method fails to converge & break & \cite{dai2002,mascarenhas2004} \\
Adam & a convex problem on which the iterates do not converge to the solution & break & \cite{reddi2018} \\
Nelder--Mead & strictly convex functions of two variables on which the simplices converge to a nonstationary point & break & \cite{mckinnon1998} \\
Multi-block ADMM & a convex problem on which the direct extension to more than two blocks diverges & break & \cite{chen2016} \\
Frank--Wolfe & problems whose solution lies on the boundary of the feasible set, where convergence is sublinear & slow & \cite{lacostejulien2015} \\
Sequential quadratic programming & the Maratos effect: near the solution the merit function rejects the full step & slow & \cite{nocedalwright2006} \\
Policy iteration & Markov decision processes that need exponentially many iterations & long path & \cite{fearnley2010} \\
$k$-means (Lloyd) & point sets that need exponentially many iterations, even in the plane & long path & \cite{arthurvassilvitskii2006,vattani2011} \\
2-opt & traveling-salesman instances with exponentially many improving exchanges & long path & \cite{englert2014} \\
\bottomrule
\end{tabular}
\end{table}

\section{A dimension classification of the family}\label{sec:theorem}
The computations of \cref{sec:members} can be compared along one clean axis: how the blow-up scales
with the dimension. \Cref{tab:compare} collects the three costs of every member.

\begin{table}[t]
\renewcommand{\arraystretch}{1.3}\footnotesize\centering
\caption{The seven members quantified in general dimension $n$: the control parameter, the baseline
cost $\Nz$ of the rule on the normal instance, the guided cost $\Ng$ and the rule cost $\Nr$ on the
\Uex, and the blow-up $\Bl=\Nr/\Ng$ with its dimension class (\cref{thm:classification}). Here
$L=\ln(1/\eps)$. For the metaheuristic the costs are expected numbers of members.}
\label{tab:compare}
{\setlength{\tabcolsep}{3.6pt}\begin{tabular}{@{}llccccl@{}}
\toprule
\textbf{Member} & \textbf{Parameter} & $\Nz$ & $\Ng$ & $\Nr$ & $\Bl$ & \textbf{Class} \\
\midrule
Newton  & flatness $2m$ & $1$ & $1$ & $\approx(2m{-}1)L$ & $=\Nr$ & dimension-free \\
Gradient & roughness $\rho$ & $1$ & $1$ & $\Theta(\rho L)$ & $=\Nr$ & dimension-free \\
Simplex & $n$, scaling & $1$ & $1$ & $2^{n}-1$ & $=\Nr$ & exponential \\
Interior-point & redundancy, $n$ & $O(\sqrt n\,L)$ & $O(\sqrt n\,L)$ & $\ge2^{n}-1$ & $2^{\Theta(n)}$ & exponential \\
Branch-and-bound & parity of $n$ & $1$ & $3$ & $\ge2^{(n+1)/2}-1$ & $2^{\Theta(n)}$ & exponential \\
Metaheuristic & basins $k$ & $1$ & $1$ & $k^n$ & $=\Nr$ & exponential \\
Active-set & degeneracy & finite & finite & $\infty$ & $\infty$ & break \\
\bottomrule
\end{tabular}}
\end{table}

\begin{theorem}\label{thm:classification}\dname{Dimension classification}
For each member $A$ let $\psi_A(n)$ denote the growth order in $n$ of the blow-up $\Bl_A$ of its
\Uex, with the accuracy $\eps$ and every control parameter other than the dimension held fixed.
Then
\[
\psi_A(n)=
\begin{cases}
\Theta(1), & A\in\{\text{Newton},\ \text{gradient}\},\\[2pt]
2^{\Theta(n)}, & A\in\{\text{simplex, interior-point, branch-and-bound, metaheuristic}\},\\[2pt]
\infty, & A=\text{active-set}.
\end{cases}
\]
Consequently the family partitions into a \emph{dimension-free} class
$\{\mathrm{Newton},\mathrm{gradient}\}$ and an \emph{exponential} class
$\{\mathrm{simplex},\mathrm{interior\mbox{-}point},\mathrm{branch\mbox{-}and\mbox{-}bound},
\mathrm{metaheuristic}\}$, with active-set as the break at the boundary.
\end{theorem}
\begin{proof}
The growth orders are read from \cref{sec:members}. The Newton and gradient instances are
separable and the rule acts on each coordinate in the same way, so the blow-up,
$\approx(2m-1)\ln(\Phi(x_0)/\eps)$ and $\Theta(\rho\ln(1/\eps))$ respectively, does not depend
on~$n$. The simplex blow-up is $2^n-1$; the interior-point blow-up is at least
$(2^n-1)/O(\sqrt n\log(1/\eps))$ and at most $O(2^nn^{5/2})$ \cite{deza2008}; the branch-and-bound
blow-up is at least $(2^{(n+1)/2}-1)/3$ and at most $2^{n+1}/3$, the size of the full binary tree;
and the expected metaheuristic blow-up is $k^n$ with $k\ge2$. All four are $2^{\Theta(n)}$. The
active-set rule cost is infinite in every dimension.
\end{proof}

The theorem classifies these seven instances; it is not a statement about every worst case of the
seven methods. Two remarks place it in context.

\begin{corollary}\label{cor:localglobal}
In the dimension-free class the blow-up is governed by a \emph{local} feature of the instance
(the order of flatness, the roughness) that acts on each coordinate separately and does not
compound. In the exponential class it is governed by a \emph{global, combinatorial} feature (a count
of vertices, of tree nodes, or of basins) that multiplies across coordinates.
\end{corollary}

\begin{remark}\label{rem:polyregime}
Between the two classes lies a polynomial regime, which none of the seven instances occupies but
which is well represented in the literature. For natural families of ill-conditioned quadratics the
condition number, and with it the cost of gradient descent, grows polynomially with $n$; and cyclic
coordinate descent can be slower than its randomized counterpart by a factor of order $n^2$
\cite{sunye2021}. The split of \cref{thm:classification} is therefore between worst cases that do
not compound across coordinates, or compound only polynomially through coupling, and those that
compound multiplicatively.
\end{remark}

\begin{figure}[tbp]\centering
\includegraphics[width=0.86\textwidth]{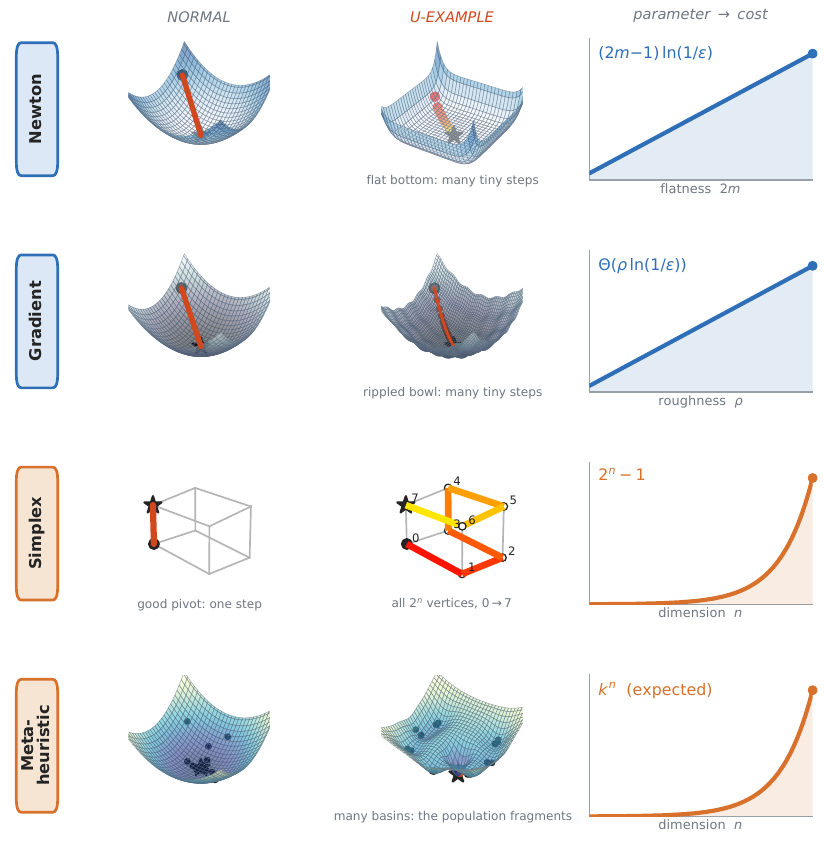}
\caption{The dimension classification, drawn out for the four members with a natural
three-dimensional geometry, in the grammar of \cref{fig:plate} (blue: dimension-free; orange:
exponential). For each, a normal instance (left) sits beside its \Uex\ (middle)---Newton's method on
a flat-bottomed bowl, gradient descent on a rippled bowl, the simplex walk on the Klee--Minty cube
(one pivot reaches the optimum; Dantzig's rule threads all $2^n$ vertices $0\to7$ to that same
vertex), and a population on a landscape with many basins. The right column plots each blow-up
against its control parameter: linear in the flatness and in the roughness, exponential in the
dimension.}
\label{fig:plate3d}
\end{figure}

\Cref{fig:plate3d} makes the dichotomy visual for the four members with a natural surface or
polytope. Each row places a three-dimensional normal instance beside its \Uex\ and, in the third
column, plots the blow-up against the member's control parameter: the dimension-free members (blue)
rise linearly with a local parameter, while the exponential members (orange) explode with the
dimension.

\section{Synthesis: one form, not one parameter}\label{sec:synthesis}
What the seven share is a \emph{form}: the optimum is near, the rule is forced the long way or
into a loop, and the cost is dominated by the detour. A natural test of how far the
unification goes is to ask for more: does a single parameter govern all seven? It does not, and the
ways in which it fails are informative.
Newton's method is governed by one integer, the flatness order, and gradient descent by one real
number, the roughness; in both the blow-up is independent of the dimension. For the simplex method
the parameter is the dimension itself, together with a scaling. For interior-point methods it is a
redundancy that must grow with the dimension. For branch-and-bound it is a switch, the parity
of~$n$; for active-set methods another switch, degeneracy, which turns a finite cost into an
infinite one. For metaheuristics it is a number of basins, and the statement is probabilistic.
The shared form, not a shared parameter, is the unifying claim.

\begin{figure}[t]\centering
\includegraphics[width=0.92\textwidth]{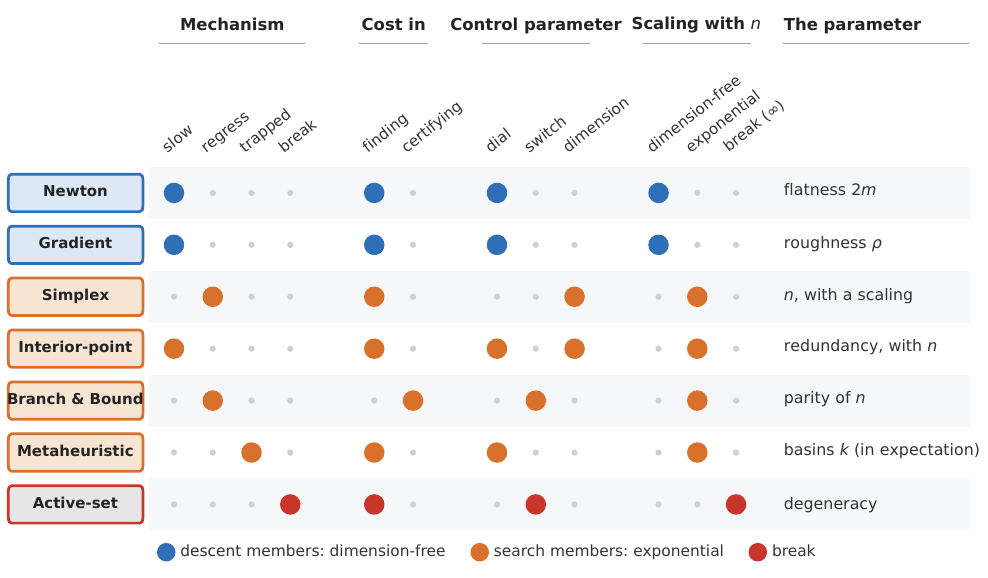}
\caption{The family along the four axes of \cref{sec:synthesis}. Each row is a member, in the order
of \cref{fig:story}; each block of columns is one axis, and a filled cell marks the value the member
takes. The descent members (blue) and the search members (orange) differ on every axis except the
shared form.}
\label{fig:synthesis}
\end{figure}

Four cross-class axes organize the family; \cref{fig:synthesis} lays the members out along them.
\begin{enumerate}
\item \textbf{The mechanism.} \Cref{prop:blowup} allows two ways of being forced the long way, and
both occur. Newton's method, gradient descent, and interior-point methods are \emph{slow}: every
step makes progress, but only a small share of what is available. The simplex method and
branch-and-bound \emph{regress}: they give back progress already made, on the cube by walking away
from the optimum, in the tree by doubling the work that remains. Metaheuristics are \emph{trapped}
with high probability, and active-set methods \emph{break}.
\item \textbf{Finding versus certifying.} For branch-and-bound the cost lies in \emph{proving} the
answer, not in reaching it. The descent members have no such separation: for them, arriving is
finishing.
\item \textbf{The kind of control parameter.} It may be a dial (flatness, roughness, number of
basins, redundancy), a switch (parity, degeneracy), or the dimension itself.
\item \textbf{The scaling with dimension.} By \cref{thm:classification} the family divides into the
forced-long \emph{descent} members (Newton, gradient: local, dimension-free) and the forced-long
\emph{search} members (simplex, interior-point, branch-and-bound, metaheuristic: global, exponential),
with active-set as the break.
\end{enumerate}
The last axis is the structural finding of the survey: the classification by dimension and the
qualitative distinction between descent and search are the same partition, reached from two
directions.

\section{Connections to other fields}\label{sec:connections}
\begin{figure}[t]\centering
\includegraphics[width=0.9\textwidth]{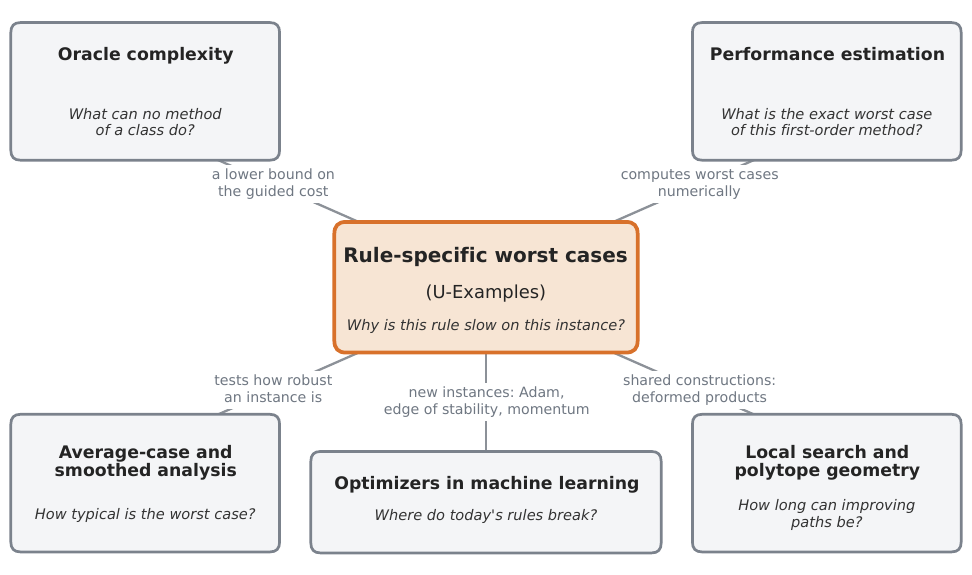}
\caption{How the rule-specific worst cases surveyed here relate to neighbouring theories. Each
theory is labelled by the question it asks of an instance or a method.}
\label{fig:connections}
\end{figure}

Rule-specific worst cases touch several neighbouring theories, each of which answers a different
question about the same instances. \Cref{fig:connections} summarizes the relations.

\subsection{Oracle complexity}
Information-based complexity asks what \emph{no} method of a class can do
\cite{nemirovskiyudin1983}. Its lower bounds are stronger than rule-specific ones in that they
quantify over all methods, and weaker in that they say nothing about why a particular good method
is slow on a particular instance. The two meet in the ill-conditioned quadratic, which is both the
classical worst case of gradient descent and, in Nesterov's chain form, the hard function for every
first-order method \cite{nesterov2004,nesterov2018}. Modern results extend the theory to nonconvex
problems \cite{carmon2020}, stochastic oracles \cite{arjevani2023}, and second-order methods
\cite{arjevani2019}. In the vocabulary of \cref{sec:framework}, an oracle lower bound is a lower
bound on the \emph{guided} cost: it shows that even the best admissible route is long. A \Uex\ is
the opposite situation, in which the guided cost is small and only the rule is at fault.

\subsection{Performance estimation}
The performance estimation problem (PEP), introduced by Drori and Teboulle \cite{drori2014} and
developed by Taylor, Hendrickx, and Glineur \cite{taylor2017,taylor2017thesis}, computes the exact
worst-case performance of a given first-order method over a class of functions as the value of a
semidefinite program, and its solution yields a worst-case function. Software makes the computation
routine \cite{taylor2017pesto,goujaud2022pepit}. PEP is the closest existing analogue of a general
theory of rule-specific worst cases, and it is constructive where the present survey is
descriptive. It differs in scope and in kind: it applies to first-order methods on convex and
related classes, it is numerical, and it targets tight rates, whereas the instances collected here
include pivoting, tree-search, and population methods, are given in closed form, and are compared
by their form. PEP-style reasoning has also settled questions about single rules. The heavy-ball
method was shown to cycle on a strongly convex function for suitable parameters by a related
computer-aided analysis \cite{lessard2016}, and it is now known that, for every choice of
parameters, heavy ball either fails to accelerate on quadratics or cycles on some smooth strongly
convex function \cite{goujaud2025}.

\subsection{Average-case and smoothed analysis}
Average-case and smoothed analysis ask how typical a worst case is. For the simplex method the
answer is that worst cases are atypical: the expected number of pivots is polynomial under natural
input distributions \cite{borgwardt1987,smale1983} and under Gaussian perturbation of arbitrary
inputs \cite{spielmanteng2004,spielmanteng2009,dadushhuiberts2020,huiberts2025,bachhuiberts2025};
see \cite{roughgarden2021} for the wider programme. For branch-and-bound the answer is mixed:
an exponential lower bound persists when Gaussian noise is added to the data
\cite{dey2023lower}. \Cref{sec:boundaries} examines how the seven instances respond to
perturbation.

\subsection{The complexity of local search and polytope geometry}
A pivot rule, a local-search heuristic, and a descent method are all procedures that follow improving
steps to a local optimum. The class PLS captures the complexity of finding such an optimum
\cite{johnson1988}, and for PLS-complete problems the standard local-search algorithm takes
exponentially many steps on some instances \cite{schafferyannakakis1991,papadimitriou1998}. The
explicit exponential instances for $k$-means \cite{arthurvassilvitskii2006,vattani2011} and for
2-opt \cite{englert2014} are \Uexs\ in all but name: good solutions are nearby, and the improvement
rule takes exponentially long to settle. For pivoting, the underlying geometric question is the
diameter of polytopes \cite{santos2012,kalaikleitman1992,todd2014}, and the constructions are
deformed products \cite{amentaziegler1999,ziegler1995}, which also link the simplex and
interior-point cubes through the curvature of the central path \cite{deza2009curvature}. For random
optimization landscapes, the overlap-gap property is a structural explanation of why local methods
stall \cite{gamarnik2021ogp}.

\subsection{Optimizers in machine learning}
The optimizers of machine learning have acquired worst cases of exactly the rule-specific kind.
Reddi, Kale, and Kumar exhibited a simple convex problem on which the Adam method fails to converge
to the optimal solution, and traced the failure to its exponential moving average \cite{reddi2018};
this is a \emph{break}, with a remedy designed against it. Gradient descent on neural networks is
observed to run at the edge of stability, with the step length in tension with the sharpest
curvature \cite{cohen2021}, the mechanism of the rippled bowl of \cref{ssec:gradient}. Momentum
methods cycle on suitable instances \cite{lessard2016,goujaud2025}, and cyclic coordinate descent has
a quadratic worst case polynomially slower than its randomized version \cite{sunye2021}. We return to
the last two in \cref{sec:outlook}; \cref{tab:beyond} lists these and further cases.

\section{From diagnosis to remedy}\label{sec:remedy}
\begin{figure}[t]\centering
\includegraphics[width=\textwidth]{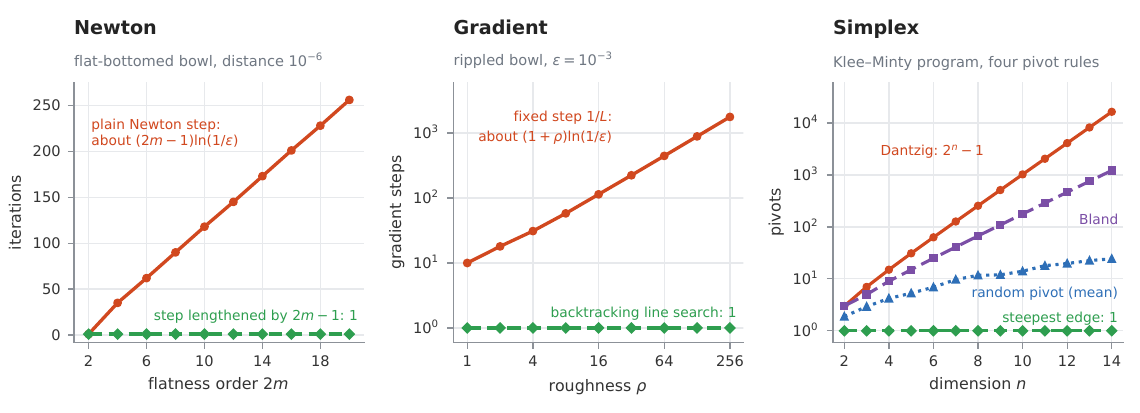}
\caption{Three remedies at work on the instances of \cref{sec:members} (computed). Each panel plots
the cost against the member's control parameter, for the rule and for a remedy. \emph{Left:} Newton
iterations to distance $10^{-6}$ on $\sum_i x_i^{2m}$: the plain step, whose cost grows linearly
with the flatness order, and the step lengthened by the factor $2m-1$, which reaches the minimizer
at once. \emph{Middle:} gradient steps on the rippled bowl \cref{eq:rippled} to accuracy $10^{-3}$:
the fixed step $1/L$, whose cost grows linearly with the roughness, and a backtracking line search
started from the unit step. \emph{Right:} pivots on the Klee--Minty program \cref{eq:km} as the
dimension grows, for Dantzig's rule, Bland's rule, a uniformly random improving pivot (mean of
$200$ runs), and the steepest-edge rule.}
\label{fig:remedy}
\end{figure}

Reading a worst case as a misdirected rule is not only diagnostic; it points at the cure. For each
member the \Uex\ names a specific way in which the rule is misled, and a standard algorithmic remedy
is the operation that removes that misdirection. \Cref{tab:remedy} is the map, and
\cref{fig:remedy} shows three of the remedies at work on the instances of \cref{sec:members}: in
each case the remedy does not merely lower the cost but removes its dependence on the control
parameter.

\begin{table}[t]
\renewcommand{\arraystretch}{1.25}\footnotesize\centering
\caption{Diagnosis to remedy. Each remedy is an operation that removes the misdirection named in the
\Uex.}
\label{tab:remedy}
\begin{tabular}{@{}p{2.2cm} >{\raggedright\arraybackslash}p{4.6cm} >{\raggedright\arraybackslash}p{5.1cm}@{}}
\toprule
\textbf{Member} & \textbf{Misdirection (the \Uex)} & \textbf{Remedy} \\
\midrule
Newton & the curvature model under-scales the step at a flat bottom & a step lengthened by the multiplicity; regularized Newton steps \\
Gradient & ripples dictate a tiny step length & line search or adaptive step lengths; smoothing; acceleration \\
Simplex & the pivot rule is lured onto a $2^n$-vertex path & a different pivot rule or a rescaling defeats \emph{this} instance; each known rule has its own \\
Interior-point & redundant walls bend the central path & presolving, which removes redundant constraints \\
Branch-and-bound & the relaxation bound prunes nothing & cutting planes; branching on general disjunctions \\
Metaheuristic & basins split the population & restarts and niching; diversity control; recombination \\
Active-set & the working set cycles at a degenerate point & anti-cycling rules; perturbation of the constraints \\
\bottomrule
\end{tabular}
\end{table}

Two patterns in \cref{tab:remedy} are worth naming. First, the remedies for the descent members
(a longer step, a line search) correct the \emph{local} information the rule reads, whereas the
remedies for the search members (cutting planes, niching, another pivot rule) change the
\emph{global} information the rule uses to prune or to diversify: the descent/search split appears
again, now in the cures. Second, the active-set remedy is of a different kind. It does not speed up
a slow process but restores termination to a broken one.

The map is not a guarantee. A remedy removes one misdirection; it does not make the method immune.
The pivot rules that defeat the Klee--Minty cube have their own exponential instances
(\cref{sec:history}), and branch-and-bound with general disjunctions still has exponential lower
bounds \cite{dey2023lower}.

\section{Boundaries and robustness of the unification}\label{sec:boundaries}
A unifying reading is useful only if its limits are stated. The \Uex\ reading organizes instances;
it is not a theorem about every worst case. This section says where the reading is exact, where it
strains, what it does not cover, and how robust the seven instances are.

\subsection{Where the reading is exact, and where it strains}\label{ssec:exact}
The reading is exact for the Newton, gradient, and simplex members: the rule is deterministic, the
potential is the natural one, and the three costs are computed in closed form. For the other four it
holds with a qualification, and each qualification is of a different kind.

\emph{Interior-point.} The blow-up is exponential in the dimension but not in the size of the input,
and the method remains polynomial (\cref{ssec:ipm}). The instance is a property of a
\emph{description} of the feasible region, not of the region itself.

\emph{Branch-and-bound.} The blow-up lies in certifying the answer, and the short route uses a move,
a general disjunction, that standard implementations do not consider. This exposes a choice that is
present for every member: the guided cost depends on the admissible set $\mathcal A$. A larger
admissible set gives a shorter guided route and a larger blow-up. We have taken $\mathcal A$ to be
the natural class of the method in each case, but the choice is a modelling decision, and it should
be read as part of the instance.

\emph{Metaheuristic.} The statement is probabilistic and holds for a model, independent descents
from random starts. Randomness means that no instance can forbid success.

\emph{Active-set.} The rule cost is infinite, and the member is a boundary case, not a ratio.

A further dependence is on the potential. Measured in function value instead of distance, the
Newton member is not a \Uex\ at all (\cref{rem:newton-f}). A \Uex\ is therefore a property of a
triple---a rule, an admissible set, and a measure of progress---and all three should be named
whenever a worst case is quoted.

\subsection{What the reading does not cover}
Not every hard instance is a \Uex, and the definition is designed to exclude two important kinds.

The first is hardness for a whole class of methods. An oracle lower bound shows that every method
of a class is slow on some function \cite{nemirovskiyudin1983,nesterov2004}; in our terms the
\emph{guided} cost is already large, so there is no misdirection to speak of. The ill-conditioned
quadratic is of this kind for the admissible set used in \cref{ssec:gradient}: steepest descent with
an exact line search, the best step along the negative gradient, still zig-zags at a rate governed
by the condition number \cite{nocedalwright2006}. This is why the gradient member is the rippled
bowl and not the narrow valley.

The second is intrinsic hardness of the problem. For several PLS-complete local-search problems
there are instances and starting solutions from which \emph{every} sequence of improving moves is
exponentially long \cite{schafferyannakakis1991}: no rule does better, and the long way is the only
way. The same is true, in a different sense, of random instances with the overlap-gap property
\cite{gamarnik2021ogp}. Such instances are hard, but they are not misdirections.

Finally, the survey is selective in two ways. The seven classes are representative, and other
methods have rule-specific worst cases of their own (\cref{tab:beyond}). And within each class the
instance was chosen as the cleanest representative: each class has other worst cases, some of a
different dimension class (\cref{rem:polyregime}), and \cref{thm:classification} is a statement
about the seven instances, not about the seven classes.

\subsection{Robustness of the constructions}\label{ssec:robust}
\begin{figure}[t]\centering
\includegraphics[width=\textwidth]{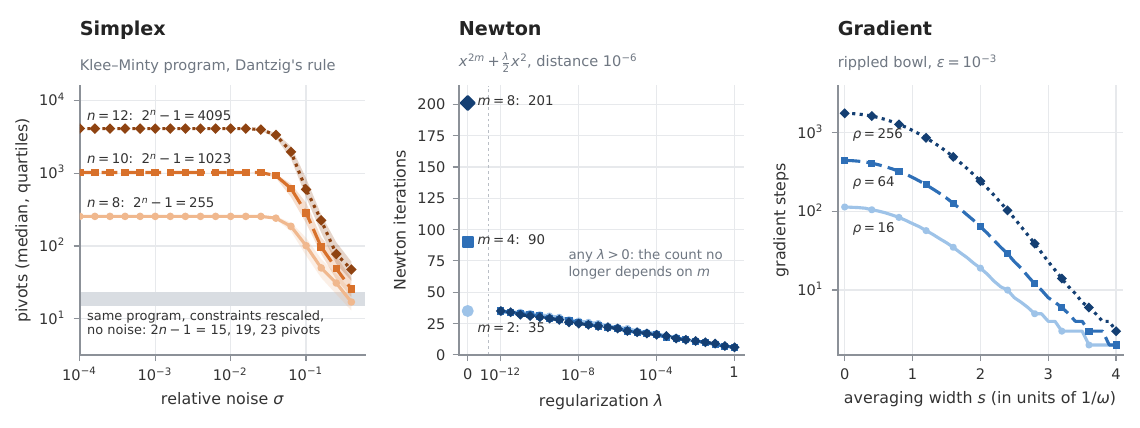}
\caption{Three members under perturbation (computed). \emph{Left:} pivots of Dantzig's rule on the
Klee--Minty program \cref{eq:km} when every nonzero coefficient is multiplied by $1+\sigma Z$ with
independent standard normal $Z$ (median and interquartile range over $60$ samples). The exponential
path survives relative noise of about $1\%$ and collapses beyond it. The grey band marks a
different experiment with no noise at all: each constraint of the same program is divided by its
norm, which changes neither the feasible region nor the optimum, and Dantzig's rule then needs only
$2n-1$ pivots ($15$, $19$, $23$).
\emph{Middle:} Newton iterations to distance $10^{-6}$ on $x^{2m}+\tfrac\lambda2x^2$ from $x_0=1$.
The unperturbed counts ($\lambda=0$, left margin) grow with $m$; any $\lambda>0$ removes the
dependence on $m$. \emph{Right:} gradient steps on the rippled bowl \cref{eq:rippled}, $\eps=10^{-3}$,
after the function is averaged over a Gaussian window of width $s$; the ripples, and the blow-up,
vanish once the window spans a few ripple wavelengths.}
\label{fig:fragility}
\end{figure}

It is often said that worst-case instances are fragile, in the sense that a small perturbation
destroys them. Smoothed analysis makes this precise for the simplex method
(\cref{sec:connections}). For the seven instances surveyed here the picture is more varied, and
\cref{fig:fragility} shows three computed examples. 

\emph{Simplex.} The Klee--Minty program is not destroyed by small noise. In the computation of
\cref{fig:fragility} (left), Dantzig's rule still takes exactly $2^n-1$ pivots for $n=8,10,12$ in
at least three quarters of the samples when every coefficient is perturbed by a relative $1\%$; the
count starts to fall near $2.5\%$ and collapses beyond: at $10\%$ the median for $n=12$ is $598$
pivots instead of $4095$, and at $25\%$ it is $78$. The instance is fragile in two other ways.
Rescaling each constraint to unit norm, which changes neither the feasible region nor the objective,
reduces the count to $2n-1$: Dantzig's rule compares reduced costs, which depend on the units in
which the constraints are written, and the scaling of \cref{eq:km} is tuned to it. Changing the
pivot rule reduces the count to one (\cref{ssec:simplex}).
The fragility is to scaling and to the choice of rule, not to small relative noise. This does not
conflict with smoothed analysis, which concerns the shadow-vertex rule and a perturbation that is
additive and measured against the largest row of the constraint matrix: the rows of \cref{eq:km}
range over many orders of magnitude, and such a perturbation overwhelms the small ones.

\emph{Newton.} The flat bottom is fragile to noise in the strictest sense. Adding
$\tfrac\lambda2\|x\|^2$ makes the Hessian nonsingular at the minimizer, and the iteration then
crawls only down to the scale at which the quadratic term takes over, after which it converges
quadratically. In \cref{fig:fragility} (middle) the unperturbed method needs $35$, $90$, and $201$
iterations for $m=2,4,8$; with $\lambda=10^{-12}$ it needs about $35$ for all three, and with
$\lambda=10^{-6}$ about $21$. The count is then of order $\ln(1/\lambda)$ whatever the flatness.

\emph{Gradient.} The rippled bowl is robust to perturbation of its parameters: nearby values of
$\rho$ and $\omega$ give the same behaviour. But it lives at a single length scale. Averaging the
function over a Gaussian window of width $s$ multiplies the ripple amplitude by
$e^{-(s\omega)^2/2}$, and a window of three to four times $1/\omega$, far below the accuracy sought,
removes the blow-up (\cref{fig:fragility}, right).

\emph{The other members.} The interior-point instance depends on redundant constraints, and
presolving removes them. Degeneracy, the switch of the active-set member, is destroyed by an
arbitrarily small perturbation of the constraints, which is how the classical perturbation methods
prevent cycling; structured problems are nevertheless often degenerate by construction. Jeroslow's
instance rests on an exact parity of integer data, which has no small perturbations; but
exponential lower bounds for branch-and-bound are not in general fragile, since one such bound
persists when Gaussian noise is added to the data \cite{dey2023lower}. The basins of the
metaheuristic member are robust: a small perturbation of a multimodal landscape leaves its basins
in place.

\Cref{tab:robust} summarizes. Fragility is common in the family but not universal, and it comes in
different kinds: to noise (Newton, active-set), to scaling and to the choice of rule (simplex), to
the description of the problem (interior-point), and to averaging (gradient). This goes some way
toward explaining a familiar asymmetry. Flat bottoms, exactly degenerate points, and Klee--Minty
scalings are rarely met in computation, whereas weak relaxations and multimodal landscapes are
everyday difficulties of integer programming and global optimization: the first group is fragile,
the second is not.

\begin{table}[t]
\renewcommand{\arraystretch}{1.25}\footnotesize\centering
\caption{Robustness of the seven instances. ``Computed'' refers to \cref{fig:fragility}; the other
entries follow from the constructions or from the literature cited in the text.}
\label{tab:robust}
\begin{tabular}{@{}p{2.2cm} >{\raggedright\arraybackslash}p{5.0cm} >{\raggedright\arraybackslash}p{4.6cm}@{}}
\toprule
\textbf{Member} & \textbf{Under small perturbation of the data} & \textbf{What removes the worst case} \\
\midrule
Newton & fragile: any quadratic term cuts the slow phase short (computed) & regularization; a lengthened step \\
Gradient & robust in $\rho$ and $\omega$ & averaging at the ripple scale (computed); line search \\
Simplex & robust to relative noise up to about $1\%$ (computed) & rescaling; another pivot rule \\
Interior-point & depends on the description, not on the region & presolving \\
Branch-and-bound & parity is exact; other lower bounds persist under noise \cite{dey2023lower} & general disjunctions; cutting planes \\
Metaheuristic & robust: basins persist & restarts; niching; recombination \\
Active-set & fragile: degeneracy is non-generic & perturbation; anti-cycling rules \\
\bottomrule
\end{tabular}
\end{table}

\section{Open challenges and outlook}\label{sec:outlook}
The seven members were assembled one at a time, each from its own literature. The natural next
question runs the other way: given only an update rule, can one tell in advance whether it admits a
\Uex, and of which class? We do not answer this question. This section records how far the framework
reaches toward an answer, states precisely what is proved, poses the rest as a conjecture, and closes
with further open problems.

\subsection{The converse question}
Write a rule as the arg-max of a \emph{local score} $\sigma_s\colon\mathcal A(s)\to\R$, so that
$T(s)=\operatorname*{arg\,max}_{s'\in\mathcal A(s)}\sigma_s(s')$; this covers Dantzig's pivot rule,
the gradient and Newton steps, and best-bound node selection. The \emph{$\Phi$-greedy} move is the
admissible step that removes the most potential. Misdirection is the disagreement between the score
and $\Phi$-progress, which suggests measuring how far a rule can disagree.

\begin{definition}\label{def:aligned}\dname{Alignment and uniform progress}
Fix an instance $I$ and an accuracy $\eps$, and call a state \emph{active} if the rule visits it
before reaching $\Phi\le\eps$. A rule $T$ is \emph{$c$-aligned} on $I$ (for $c\ge1$) if at every
active state $s$,
$\Phi(s)-\Phi(T(s))\ge\tfrac1c\max_{s'\in\mathcal A(s)}(\Phi(s)-\Phi(s'))$; it is
\emph{$\Phi$-greedy} if it is $1$-aligned. The instance has \emph{uniform progress} $\mu\in(0,1]$ if
at every active state an admissible step can remove a fraction $\mu$ of the remaining potential,
$\max_{s'\in\mathcal A(s)}(\Phi(s)-\Phi(s'))\ge\mu\,\Phi(s)$.
\end{definition}

\begin{proposition}\label{prop:aligned}\dname{Aligned rules on uniformly decreasing potentials}
If $I$ has uniform progress $\mu$ and $T$ is $c$-aligned on $I$, then $\delta_k\le1-1/c$ at every
step and
\[
  \Nr(I;\eps)\le\Bigl\lceil\tfrac{c}{\mu}\,\ln\tfrac{\Phi(s_0)}{\eps}\Bigr\rceil ,
\]
so $\Bl(I;\eps)$ obeys the same bound, independently of the dimension and of every other
feature of the instance.
\end{proposition}

The proof is in \cref{app:proofs}. The proposition turns the misdirection of \cref{def:misdirection}
into a statement about rules: a \Uex\ needs either an unbounded misalignment $c$ or a vanishing
uniform progress $\mu$. The two descent members are of the first kind, and the proposition
reproduces their costs. On the flat bottom a lengthened Newton step reaches the minimizer
($\mu=1$), while the plain step removes only the fraction $1/(2m-1)$ of the distance, so $c=2m-1$
and the bound is $(2m-1)\ln(\Phi(x_0)/\eps)$. On the rippled bowl, along the trajectory from the
symmetric start, a suitable step along the negative gradient again reaches the minimizer
($\mu=1$), and $c\le2(1+\rho)$, which gives the upper bound of \cref{lem:rippled}. The Klee--Minty cube lies
outside the proposition altogether, because its pivots regress ($c=\infty$), in line with the
detour mechanism of \cref{prop:blowup}.

\subsection{Gadgets, composition, and a conjecture}
Each member appears to be a single low-dimensional configuration repeated across coordinates. We
make this reading explicit.

\begin{definition}\label{def:gadget}\dname{Misdirection gadget}
A \emph{gadget} for $T$ is a finite instance $G$ with $\Nr(G)>\Ng(G)$, on which the rule is
misdirected on some step. A gadget is \emph{cycling} if the rule never reaches the target,
$\Nr(G)=\infty$.
\end{definition}

\begin{definition}\label{def:composition}\dname{Additive and coupled composition}
The \emph{additive product} $I\otimes J$ has state space $S_I\times S_J$, potential $\Phi_I+\Phi_J$,
one-coordinate admissible moves, and target $\Sstar_I\times\Sstar_J$, so that the coordinates do not
interact. A \emph{coupled product} is the same, except that the admissible moves in one coordinate
depend on the state of the others---a deformed product in the sense of \cite{amentaziegler1999}.
\end{definition}

When a rule acts coordinatewise, an additive product of a fixed gadget has the same blow-up as the
gadget itself: the guided route and the rule both pay once per coordinate. Dimension then leaves the
blow-up untouched, which is the behaviour of the dimension-free class of
\cref{thm:classification}, where the blow-up is set by a local parameter such as the flatness order.
A coupled product can do more: in the Klee--Minty cube the coupling lets a rule be lured through all
$2^n$ joint configurations. These observations suggest the following statement, which we have not
been able to prove.

\begin{conjecture}\label{conj:general}
Let $T$ be a rule whose score depends on the state only through a finite local description, the same
at every state. Then $T$ admits a \Uex\ family if and only if it admits a gadget that can be composed
across coordinates. Moreover, additive composition yields a blow-up that does not compound with the
dimension, coupled composition can yield a blow-up exponential in the dimension, and a cycling gadget
yields a break, so that the trichotomy of \cref{thm:classification} would be the trichotomy of the
rule's gadgets.
\end{conjecture}

Two gaps separate the conjecture from a theorem. First, coupling is not sufficient on its own: which
couplings force exponentially long rule paths, rather than merely permit them, is precisely the
pivot-rule-specific question that each exponential simplex construction answers by hand
(\cref{sec:history}). Second, the \emph{only if} direction would require extracting a gadget from an
arbitrary blow-up family, and we know of no general extraction argument; for adaptive or
history-dependent rules, such as Zadeh's \cite{disser2023zadeh}, even the statement needs
reformulating.

\subsection{A heuristic reading of a rule}
\begin{figure}[t]\centering
\includegraphics[width=\textwidth]{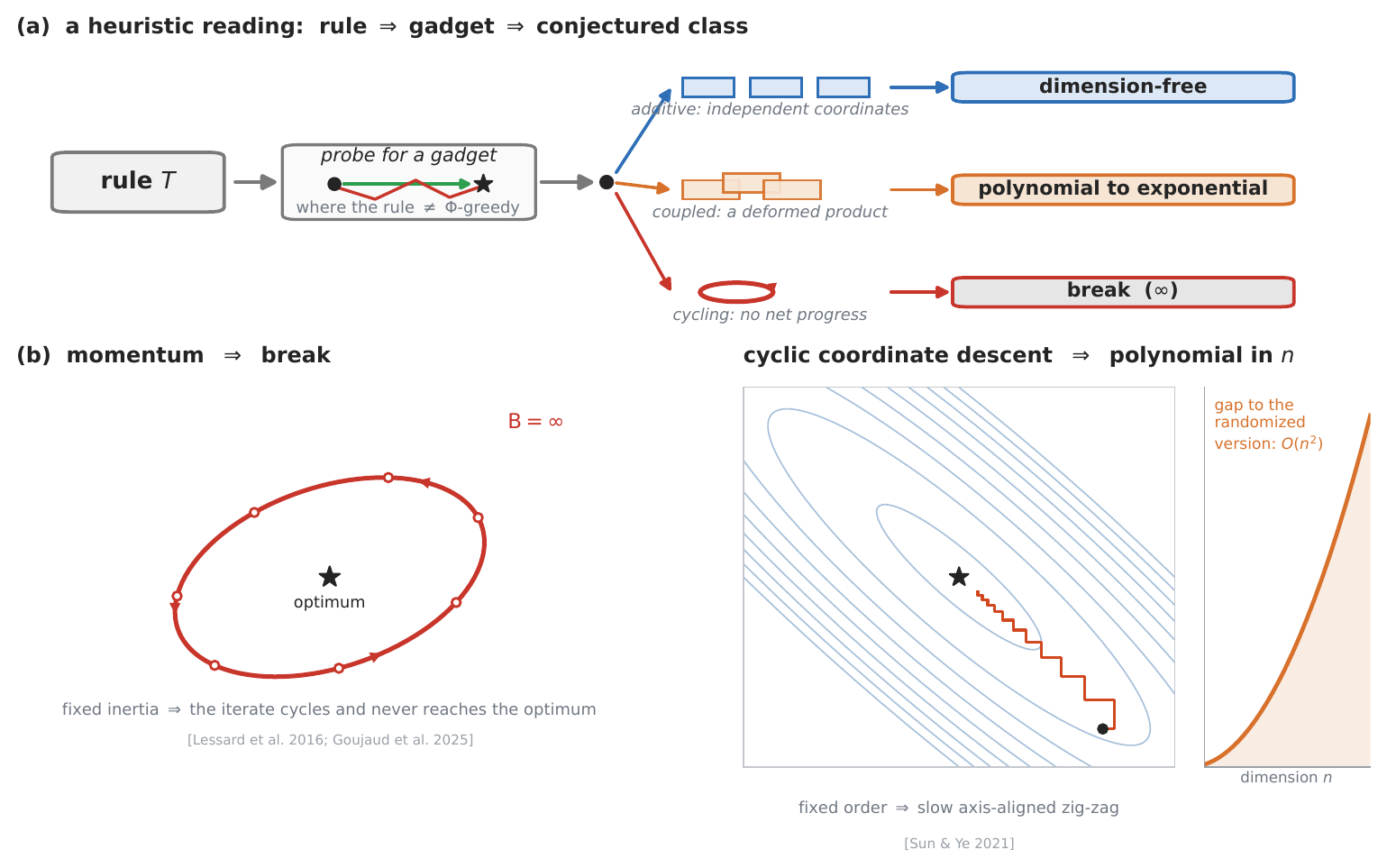}
\caption{The heuristic reading of \cref{sec:outlook}. \emph{(a)} A rule is probed for a misdirection
gadget, and the way the gadget composes---additively (independent tiles), coupled (a deformed
product), or cycling---suggests a class: dimension-free, polynomial to exponential, or break
(\cref{conj:general}). \emph{(b)} Two methods outside the seven, read this way: heavy-ball momentum
settles into a cycle that never reaches the optimum (a break, $\Bl=\infty$;
\cite{lessard2016,goujaud2025}), and cyclic coordinate descent zig-zags slowly through an
ill-conditioned valley, with a gap to its randomized version that is polynomial in $n$
\cite{sunye2021}.}
\label{fig:constructor}
\end{figure}

The conjecture suggests a three-step heuristic: \emph{(i)} probe the rule for a state where its
score's arg-max differs from the $\Phi$-greedy move; \emph{(ii)} ask whether the resulting gadget
composes additively, couples, or cycles; \emph{(iii)} read off the expected class.
\Cref{fig:constructor}(a) draws the pipeline. Run on one member from each regime
(\cref{tab:consistency}), it lands on the classes of \cref{thm:classification}. This is a consistency
check rather than evidence for the conjecture, since the seven members informed the reading.

\begin{table}[htbp]
\renewcommand{\arraystretch}{1.25}\centering\footnotesize
\caption{The heuristic reading run on one member from each regime. The probe locates where the rule
leaves $\Phi$-greedy, and how the resulting gadget composes suggests the class. The seven members
informed the reading, so this is a consistency check, not a test.}
\label{tab:consistency}
\begin{tabular}{@{}p{1.7cm} >{\raggedright\arraybackslash}p{4.6cm} >{\raggedright\arraybackslash}p{3.2cm} >{\raggedright\arraybackslash}p{2.0cm}@{}}
\toprule
\textbf{Member} & \textbf{Probe: where the rule leaves $\Phi$-greedy} & \textbf{Composition (why)} & \textbf{Predicted class} \\
\midrule
Newton & the curvature-scaled step undershoots at a flat minimum & additive: coordinates separate & dimension-free \\
Simplex & Dantzig's rule enters the variable with the largest coefficient, away from the optimum & coupled: a deformed product & exponential \\
Active-set & the rule revisits a working set at a degenerate point & cycling: no net progress & break \\
\bottomrule
\end{tabular}
\end{table}

\subsection{Two methods outside the seven}
Two methods that played no part in assembling the seven are consistent with the reading
(\cref{fig:constructor}(b)).

\paragraph{Momentum (heavy ball)}
For Polyak's heavy-ball method at fixed parameters, the inertial term can overshoot the
$\Phi$-greedy step, and Lessard, Recht, and Packard exhibited a one-dimensional strongly convex
function on which the iterates cycle rather than converge \cite{lessard2016}. Goujaud, Taylor, and
Dieuleveut showed that this is no accident: for every condition number and every choice of
parameters, heavy ball is either not accelerated on quadratics or admits a smooth strongly convex
function on which it cycles \cite{goujaud2025}. The gadget cycles, and heavy ball sits with active-set
at the break boundary of the family.

\paragraph{Cyclic coordinate descent}
For cyclic coordinate descent the natural gadget is the update order: a coordinate the fixed order
reaches last, so that progress on it is partly undone by the intervening updates through the
off-diagonal coupling of a quadratic. Sun and Ye constructed convex quadratics on which cyclic
coordinate descent is $O(n^2)$ times slower than its randomized counterpart, for any fixed update
order \cite{sunye2021}. The gadget is coupled, and here the coupling compounds only polynomially:
the slowdown is polynomial in the dimension, the regime of \cref{rem:polyregime}.

\subsection{Further open problems}
Beyond \cref{conj:general}, the survey leaves several questions open.
\begin{enumerate}
\item \emph{Adaptive and stochastic rules.} The framework treats randomized methods only through
expectations, and adaptive methods---quasi-Newton updates, step-length adaptation, the optimizers of
machine learning---change their local description along the trajectory. They have worst cases of
their own, such as the non-convergence of Adam \cite{reddi2018}, but a misdirection theory for such
rules is missing.
\item \emph{A theory of robustness.} Smoothed analysis explains why the simplex member is rarely
met \cite{spielmanteng2004}, whereas an exponential lower bound for branch-and-bound survives
random perturbation \cite{dey2023lower}. \Cref{sec:boundaries} shows that the seven instances differ widely
in how fragile they are. When is a \Uex\ fragile, and to what?
\item \emph{Rule-independent lower bounds.} The exponential bound of \cite{dissermosis2025} holds for
every pivot rule of the active-set method, at the price of a polynomial objective. Whether the
simplex method, with a linear objective, admits a polynomial pivot rule remains the central open
problem of the area.
\item \emph{Tightness.} For several members the \Uex\ lower bound and the best known upper bound for
the same rule differ by polynomial factors; closing these gaps would sharpen
\cref{thm:classification}.
\item \emph{Computing worst cases.} Performance estimation computes exact worst cases for first-order
methods \cite{drori2014,taylor2017}. An analogue for discrete rules---pivoting, branching, local
search---would turn the catalogue of this survey into a procedure.
\end{enumerate}

\section{Conclusion}\label{sec:conclusion}
The worst-case instances of optimization were built one algorithm class at a time, over seventy
years and in separate literatures. Surveying them together, we found that the rule-specific ones share a common form: an optimum placed near the start, and an update rule misdirected onto a long trajectory, or a loop, to reach it. We called such an instance a \Uex. This reading, and the framework that makes it precise, are the main contributions of the survey. The framework treats every algorithm as a rule on
a state space with a progress potential, so that misdirection and blow-up are defined once for all
seven classes, and a simple accounting bound names the two ways in which a rule can be forced the
long way: by wasting the progress on offer, or by giving back progress already made. Computed in
general dimension, the blow-up is independent of the dimension for the two descent members and
exponential for the four search members, with degenerate cycling as a breakdown; this partition
coincides with the qualitative distinction between descent and search.

The practical message is threefold. First, a worst case is a diagnosis: it names the way a rule
is misdirected, and the standard remedy is the operation that removes it (\cref{sec:remedy}).
Second, a worst case should be quoted together with its rule, its admissible moves, and its
measure of progress, because changing any of the three can make it disappear
(\cref{ssec:exact}). Third, worst cases are not equally fragile: flat bottoms, exact
degeneracy, and Klee--Minty scalings vanish under small changes, whereas weak relaxations and
multimodal landscapes do not (\cref{ssec:robust}). Conceptually, the comparison across classes isolates what is not shared: there is no common
control parameter, only a common form.

The reading also has limits, which we have tried to state plainly. It is exact for deterministic
rules with a clear potential, probabilistic for population methods, and degenerate for cycling.
Several instances are classical and are used with attribution; the framework, the classification,
and the cross-class comparison are what this survey adds to them. The main open question is the
converse one: whether the update rule alone determines the existence and the class of a worst case.
We recorded what can be proved, an upper bound for well-aligned rules, stated the rest as a
conjecture, and listed further problems that we hope will attract readers from each of the
communities whose constructions are collected here.

\appendix
\section{Proofs}\label{app:proofs}

\begin{proof}[Proof of \cref{prop:blowup}]
Write $\Delta^{\mathrm r}_k=\Phi(s_k)-\Phi(s_{k+1})$ for the progress of the $k$th rule step. Since
the rule's choice is admissible, $\Delta^{\mathrm r}_k\le\Delta^{\mathrm g}_k$, and by
\cref{def:misdirection}, $\max\{0,\Delta^{\mathrm r}_k\}=(1-\delta_k)\Delta^{\mathrm g}_k\le
(1-\delta_k)\Delta^\ast$. Splitting each step into its gain and its regress,
\[
  \Phi(s_0)-\Phi(s_{\Nr})
  =\sum_{k<\Nr}\max\{0,\Delta^{\mathrm r}_k\}-R
  \le\Delta^\ast\sum_{k<\Nr}(1-\delta_k)-R
  =\Delta^\ast\,\Nr\,(1-\bar\delta)-R .
\]
Because $\Phi(s_{\Nr})\le\eps$, rearranging gives the bound on $\Nr$. For the blow-up, every
admissible step of the guided route removes at most $\Delta^\ast$, so
$\Ng\ge(\Phi(s_0)-\eps)/\Delta^\ast$, that is, $\kappa\ge1$; dividing the bound on $\Nr$ by
$\Ng=\kappa(\Phi(s_0)-\eps)/\Delta^\ast$ gives the bound on $\Bl$.
\end{proof}

\begin{lemma}\label{lem:km}\dname{The deformed cube}
Let $n\ge1$ and $\tau\in(0,\tfrac12)$, and consider the polytope
$\{x\in\R^n:\ 0\le x_1\le1,\ \tau x_{j-1}\le x_j\le1-\tau x_{j-1}\ (j=2,\dots,n)\}$ with objective
$\max x_n$. Label its vertices by $b\in\{0,1\}^n$, where $b_j=1$ when the upper constraint on $x_j$
is tight. Then:
\emph{(a)} the unique maximizer of $x_n$ is the vertex $e=(0,\dots,0,1)$, which is adjacent to the start
$0^n$;
\emph{(b)} every vertex $b\ne e$ has an improving neighbour whose label is one Hamming step closer
to $e$;
\emph{(c)} the vertices in reflected Gray-code order, $G_1=(0,1)$ and
$G_n=(G_{n-1}\times\{0\},\ \overline{G_{n-1}}\times\{1\})$ with $\overline{G}$ the reversed list, form
a path along which $x_n$ strictly increases, from $0^n$ to $e$ through all $2^n$ vertices.
\end{lemma}

\begin{proof}
At a vertex, $x_j=\tau x_{j-1}$ if $b_j=0$ and $x_j=1-\tau x_{j-1}$ if $b_j=1$ (with $x_0=0$), so a
change in $x_j$ alters $x_{j+1}$ by a factor $+\tau$ or $-\tau$ according as $b_{j+1}=0$ or $1$.
Flipping $b_j$ from $0$ to $1$ increases $x_j$, and therefore increases $x_n$ exactly when
$b_{j+1},\dots,b_n$ contain an even number of ones; flipping $b_j$ from $1$ to $0$ does the opposite.
Labels differing in one coordinate are adjacent vertices, since for $\tau<\tfrac12$ the polytope is
combinatorially a cube.
\emph{(a)} $x_n\le1-\tau x_{n-1}\le1$, with equality only if $x_{n-1}=0$, which forces
$b_{n-1}=0$ and $x_{n-2}=0$, and recursively $b_1=\dots=b_{n-1}=0$; the label $e$ differs from $0^n$
in one coordinate.
\emph{(b)} If $b_n=0$, flipping $b_n$ to $1$ is improving (no ones follow position $n$) and moves one
step closer to $e$. If $b_n=1$ and $b\ne e$, let $j<n$ be the largest index with $b_j=1$; exactly one
entry equal to $1$ follows position $j$, so flipping $b_j$ to $0$ is improving, and it moves one step closer to $e$.
\emph{(c)} By induction on $n$. On the first half of $G_n$, $b_n=0$ and $x_n=\tau x_{n-1}$, which
increases along $G_{n-1}$ by the induction hypothesis applied to $x_{n-1}$; the middle step flips
$b_n$ from $0$ to $1$, which is improving; on the second half, $b_n=1$ and $x_n=1-\tau x_{n-1}$,
which increases because $x_{n-1}$ decreases along the reversed list. The list starts at $0^n$, ends at
$e$, and contains every label once.
\end{proof}

\begin{lemma}\label{lem:rippled}\dname{The rippled bowl}
Let $f$ be the function \cref{eq:rippled} with $\rho\ge1$ and $\omega>0$, let
$2\rho/\omega\le\eps<1$, and let $x_0=(1,\dots,1)$ and $\Phi(x)=\|x\|_\infty$. Then:
\emph{(a)} $\nabla f$ is Lipschitz with constant $L=1+\rho$;
\emph{(b)} the unit step $x_0-\nabla f(x_0)$ has potential at most $\rho/\omega\le\eps/2$;
\emph{(c)} gradient descent with step length $1/L$ from $x_0$ first reaches $\Phi\le\eps$ after
$\Nr$ steps, where
\[
  \tfrac{2\rho-1}{3}\,\ln\tfrac1\eps\ \le\ \Nr\ \le\ \bigl\lceil 2(1+\rho)\ln\tfrac1\eps\bigr\rceil .
\]
\end{lemma}

\begin{proof}
The function is separable, with partial derivatives
$g(x_i)=x_i-\tfrac\rho\omega\sin(\omega x_i)$ and second derivatives
$1-\rho\cos(\omega x_i)\in[1-\rho,\,1+\rho]$. Since $\rho\ge1$, the second derivatives are bounded
in absolute value by $1+\rho$, and the bound is attained; this is~(a). For~(b), each coordinate of
the unit step is $x_i-g(x_i)=\tfrac\rho\omega\sin(\omega x_i)$, of absolute value at most
$\rho/\omega$.

For~(c), all coordinates of $x_0$ are equal and the iteration acts on each in the same way, so it
is enough to follow one coordinate, $t_{k+1}=t_k-\eta\,g(t_k)$ with $\eta=1/(1+\rho)\le\tfrac12$ and
$t_0=1$, and $\Phi(x_k)=|t_k|$. Suppose $|t|\ge2\rho/\omega$. Then
$\tfrac\rho\omega|\sin(\omega t)|\le\tfrac12|t|$, so $g(t)$ has the sign of $t$ and
$\tfrac12|t|\le|g(t)|\le\tfrac32|t|$. Because $\tfrac32\eta\le\tfrac34$, the next iterate has the
sign of $t$ and
\[
  \bigl(1-\tfrac32\eta\bigr)|t|\ \le\ |t-\eta\,g(t)|\ \le\ \bigl(1-\tfrac12\eta\bigr)|t| .
\]
Both inequalities apply at every step $k<\Nr$, because $|t_k|>\eps\ge2\rho/\omega$ there. The upper
one gives $|t_k|\le(1-\tfrac12\eta)^k\le e^{-k\eta/2}$, which is at most $\eps$ once
$k\ge2(1+\rho)\ln(1/\eps)$. The lower one gives
$\eps\ge|t_{\Nr}|\ge(1-\tfrac32\eta)^{\Nr}$, hence
$\Nr\ge\ln(1/\eps)\big/\bigl(-\ln(1-\tfrac32\eta)\bigr)$; and $-\ln(1-u)\le u/(1-u)$ with
$u=\tfrac32\eta$ gives $-\ln(1-\tfrac32\eta)\le3/(2\rho-1)$.
\end{proof}

\begin{lemma}\label{lem:jeroslow}\dname{Jeroslow's integer program}
Let $P_n=\{x\in[0,1]^n:\ 2\sum_ix_i=n\}$.
\emph{(a)} If $n$ is odd, $P_n$ contains no integer point.
\emph{(b)} If $n$ is odd and at most $(n-1)/2$ of the variables are fixed at values in $\{0,1\}$,
the resulting face of $P_n$ is nonempty. Consequently branch-and-bound with variable branching on
\cref{eq:jeroslow} processes at least $2^{(n+1)/2}-1$ nodes.
\emph{(c)} If $n$ is odd, the two sides of the disjunction $\sum_ix_i\le(n-1)/2$ or
$\sum_ix_i\ge(n+1)/2$ are both empty in $P_n$, so branching on it at the root closes the tree with
three nodes.
\emph{(d)} If $n$ is even, every vertex of $P_n$ is integral.
\end{lemma}

\begin{proof}
\emph{(a)} For integer $x$ the left-hand side $2\sum_ix_i$ is even, and $n$ is odd.
\emph{(b)} Fix $a$ variables at $1$ and $z$ at $0$, with $a+z=d\le(n-1)/2$. The remaining $n-d$
variables must sum to $n/2-a$ with each in $[0,1]$, which is possible exactly when
$0\le n/2-a\le n-d$, that is, when $a\le n/2$ and $z\le n/2$; both hold. The relaxation at every
node of depth $d\le(n-1)/2$ is therefore feasible, and by~(a) it has no integral solution. No
incumbent is ever found, so no node is discarded by its bound. Every node of depth at most
$(n-1)/2$ is thus processed and branched, and there are $\sum_{d=0}^{(n-1)/2}2^d=2^{(n+1)/2}-1$ of
them.
\emph{(c)} Every integer point satisfies one side of the disjunction, and both sides contradict
$\sum_ix_i=n/2$.
\emph{(d)} A vertex of $P_n$ has at most one coordinate strictly between $0$ and $1$, since $P_n$ is
cut from the cube by a single equation; the other coordinates are integers and sum, with it, to the
integer $n/2$, so it is an integer too.
\end{proof}

\begin{proof}[Proof of \cref{prop:aligned}]
At a visited state $s_k$ with $\Phi(s_k)>\eps$, alignment and uniform progress give
$\Phi(s_k)-\Phi(s_{k+1})\ge\Delta^{\mathrm g}_k/c\ge(\mu/c)\Phi(s_k)$, so that
$\delta_k=1-\Delta^{\mathrm r}_k/\Delta^{\mathrm g}_k\le1-1/c$ and
$\Phi(s_{k+1})\le(1-\mu/c)\,\Phi(s_k)$. By induction,
$\Phi(s_k)\le(1-\mu/c)^k\Phi(s_0)\le e^{-k\mu/c}\,\Phi(s_0)$, which is at most $\eps$ once
$k\ge(c/\mu)\ln(\Phi(s_0)/\eps)$. Hence $\Nr\le\lceil(c/\mu)\ln(\Phi(s_0)/\eps)\rceil$, and since
$\Ng\ge1$ the same number bounds $\Bl$.
\end{proof}

\section*{Declaration on the use of AI tools}
During the preparation of this manuscript, the author used a generative AI tool to assist with language editing, refinement, and improvement of the Introduction. The scientific content, including the original idea,  architecture, experiments, results and conclusions, was developed by the author. All AI-assisted material was reviewed, revised, and approved by the author. The author assumes responsibility for all content, including its accuracy, originality, and integrity.

\bibliographystyle{siamplain}
\bibliography{references}
\end{document}